\documentclass[twoside,a4paper,reqno,11pt]{amsart} 
\usepackage{amsfonts, amsbsy, amsmath, amssymb, latexsym}
\usepackage{mathrsfs,array}
\usepackage[top=28mm,right=28mm,bottom=28mm,left=28mm]{geometry}
\usepackage{stmaryrd}
\usepackage{bm}

\usepackage{hyperref}

\renewcommand{\a}{\alpha}
\renewcommand{\b}{\beta}

\newcommand{\e}{\epsilon}
 
 \renewcommand{\O}{\Omega}

 \renewcommand{\to}{\rightarrow}
 \newcommand{\s}{\sigma}

\newcommand{\G}{\bar{G}}

\newcommand{\la}{\langle}
\newcommand{\ra}{\rangle}

\newcommand{\leqs}{\leqslant}
\newcommand{\geqs}{\geqslant}
 
\newcommand{\what}{\widehat} 
 \newcommand{\vs}{\vspace{3mm}}

\makeatletter
\newcommand{\imod}[1]{\allowbreak\mkern4mu({\operator@font mod}\,\,#1)}
\makeatother

\newtheorem{theorem}{Theorem} 
\newtheorem*{conj*}{Conjecture}

\newtheorem{corol}[theorem]{Corollary}

\newtheorem{thm}{Theorem}[section] 
\newtheorem{prop}[thm]{Proposition} 
\newtheorem{lem}[thm]{Lemma}

\theoremstyle{definition}
\newtheorem{rem}[thm]{Remark}
\newtheorem{remk}{Remark}

\newtheorem*{def-non}{Definition}
\newtheorem*{rem-non}{Remark}

\begin{document}

 \author{Timothy C. Burness}
 \address{School of Mathematics, University of Bristol, Bristol BS8 1TW, UK}
 \email{t.burness@bristol.ac.uk}
  
\title{On base sizes for almost simple primitive groups}

\begin{abstract}
Let $G \leqs {\rm Sym}(\Omega)$ be a finite almost simple primitive permutation group, with socle $G_0$ and point stabilizer $H$. A subset of $\O$ is a base for $G$ if its pointwise stabilizer is trivial; the base size of $G$, denoted $b(G)$, is the minimal size of a base.
We say that $G$ is standard if $G_0 = A_n$ and $\O$ is an orbit of subsets or partitions of $\{1, \ldots, n\}$, or if $G_0$ is a classical group and $\O$ is an orbit of subspaces (or pairs of subspaces) of the natural module for $G_0$. The base size of a standard group can be arbitrarily large, in general, whereas the situation for non-standard groups is rather more restricted. Indeed, we have $b(G) \leqs 7$ for every non-standard group $G$, with equality if and only if $G$ is the Mathieu group ${\rm M}_{24}$ in its natural action on $24$ points. In this paper, we extend this result by classifying the non-standard groups with $b(G)=6$. The main tools include recent work on bases for actions of simple algebraic groups, together with probabilistic methods and improved fixed point ratio estimates for exceptional groups of Lie type.
\end{abstract}


\date{\today}
\maketitle
 
\section{Introduction}\label{s:intro}

Let $G \leqs {\rm Sym}(\O)$ be a permutation group and recall that a subset of $\O$ is a base for $G$ if its pointwise stabilizer in $G$ is trivial. The base size of $G$, denoted by $b(G)$, is the minimal cardinality of a base. For example, if $V$ is a finite dimensional vector space and one views $G = {\rm GL}(V)$ as a permutation group on $V$, then a subset of $V$ is a base if and only if it contains a basis, hence $b(G)$ coincides with the dimension of $V$ in this case.  

The study of bases has a long history, particularly in the context of finite primitive permutation groups. For example, notice that if $\O$ is finite then $|G| \leqs |\O|^{b(G)}$, so an upper bound on $b(G)$ provides a bound on the order of $G$ in terms of its degree. The investigation of bounds of this form for finite primitive groups has been a widely studied problem for more than a century, with numerous applications (see \cite{Babai, Bochert, L10, Maroti, PS}, for example). In more recent decades, following the seminal work of Sims \cite{Sims} in the early 1970s, bases have been used extensively in the computational study of finite permutation groups (see \cite{Seress_book} for more details). Bases also arise naturally in several other settings, including graph theory and combinatorics, and we refer the reader to \cite{BCam} for further background and connections to other measures on groups and graphs.
 
In this paper, we focus on bases for finite almost simple primitive groups.  In studying the base sizes of such groups, it is natural to make a distinction between so-called \emph{standard} and \emph{non-standard} groups, according to the following definition.

\begin{def-non}\label{d:ns}
Let $G \leqs {\rm Sym}(\O)$ be a finite almost simple primitive group, with socle $G_0$ and point stabilizer $H$. We say that $G$ is \emph{standard} if one of the following holds:
\begin{itemize}\addtolength{\itemsep}{0.2\baselineskip}
\item[{\rm (i)}] $G_0 = A_n$ and $\O$ is an orbit of subsets or partitions of $\{1, \ldots, n\}$; or
\item[{\rm (ii)}] $G_0$ is a classical group with natural module $V$ and either $\O$ is an orbit of subspaces (or pairs of subspaces) of $V$, or $G_0={\rm Sp}_{n}(q)$, $q$ is even and $H \cap G_0 = O_{n}^{\pm}(q)$.
\end{itemize}
Otherwise, $G$ is \emph{non-standard}.
\end{def-non}

Typically, if $G$ is standard then the point stabilizer $H$ has large order, in the sense that $|G|$ is not bounded above by a polynomial in $|G:H|$ of fixed degree (for example, if $G = {\rm L}_{n}(q)$ then $H$ is a maximal parabolic subgroup of $G$). In particular, this implies that the base size of a standard group can be arbitrarily large. For example, it is easy to see that $b(G) = n+1$ for the natural action of $G={\rm PGL}_{n}(q)$ on the set of $1$-dimensional subspaces of the natural module. In stark contrast, a highly influential conjecture of Cameron and Kantor \cite{CK} from 1993 predicts the existence of an absolute constant $c$ such that $b(G) \leqs c$ for every non-standard group $G$. This conjecture was proved by Liebeck and Shalev \cite{LSh2}, revealing a remarkable dichotomy for base sizes of almost simple primitive groups. Moreover, through the later work of several authors \cite{Bur7,BGS,BLS,BOB}, it has been shown that $c=7$ is the optimal constant, with $b(G) = 7$ if and only if $G = {\rm M}_{24}$ in its natural $5$-transitive action on $24$ points (this confirms a conjecture of Cameron \cite[p.122]{CamPG}).

One of the main open problems in this area is to determine the precise base size of all non-standard permutation groups. The groups with an alternating or sporadic socle are handled in \cite{BLS, BOB} (also see \cite{NNOW}), and there are partial results for classical groups \cite{Bur7,BGS3,BGS2} and exceptional groups of Lie type \cite{BLS}. With a view towards applications, there has been a focus on the case $b(G)=2$, but a complete classification remains out of reach. More precisely, $k=7$ is the only integer in the range $\{2, \ldots, 7\}$ for which we have a complete classification of the non-standard groups $G$ with $b(G)=k$.

The aim of this paper is to make further progress on this problem. In recent joint work with Guralnick and Saxl \cite[Theorem 13]{BGS4}, we used results on bases for actions of simple algebraic groups to show that there are infinitely many non-standard groups $G$ with $b(G)=6$. By appealing to the connection with bases for algebraic groups, as well as the earlier work in \cite{Bur7,BGS,BLS,BOB} and the probabilistic methods originally introduced in \cite{LSh2}, we can classify the groups with $b(G)=6$. The following is our main result. 

\begin{theorem}\label{t:main}
Let $G \leqs {\rm Sym}(\O)$ be a finite almost simple primitive non-standard permutation group with socle $G_0$ and point stabilizer $H$. Then $b(G)=6$ if and only if one of the following holds:
\begin{itemize}\addtolength{\itemsep}{0.2\baselineskip}
\item[{\rm (i)}] $(G,H) = ({\rm M}_{23}, {\rm M}_{22})$, $({\rm Co}_{3}, {\rm McL}.2)$, $({\rm Co}_{2}, {\rm U}_{6}(2).2)$ or $({\rm Fi}_{22}.2, 2.{\rm U}_{6}(2).2)$.
\item[{\rm (ii)}] $G_0 = E_7(q)$ and $H = P_7$. 
\item[{\rm (iii)}] $G_0 = E_6(q)$ and $H = P_1$ or $P_6$.
\end{itemize}
\end{theorem}

\begin{remk} 
In Theorem \ref{t:main} we write $P_i$ to denote a standard maximal parabolic subgroup corresponding to deleting the $i$-th node in the Dynkin diagram of $G_0$ (where we label the Dynkin diagram in the usual way, following Bourbaki \cite{Bou}). In parts (ii) and (iii), we have $b(G)=6$ for all $q$ and for all appropriate almost simple groups $G$ with the given socle (in (iii), note that the maximality of $H$ implies that $G$ does not contain any graph or graph-field automorphisms of $G_0$).
\end{remk}

By combining \cite[Theorems 3 and 4]{BLS} with the probabilistic proof of Theorem \ref{t:main}, we obtain the following asymptotic result on the abundance of bases of size $5$ for almost simple exceptional groups. This is a refinement of \cite[Theorem 2]{BLS}.

\begin{corol}\label{t:main2}
Let $G_n \leqs {\rm Sym}(\O_n)$, $n \in \mathbb{N}$, be a sequence of finite almost simple primitive groups with socle an exceptional group of Lie type. Assume that $|G_n| \to \infty$ as $n \to \infty$. Also assume that the sequence contains only finitely many groups of the form given in parts (ii) and (iii) of Theorem $\ref{t:main}$. Then the probability that $5$ randomly chosen points in $\O_n$ form a base for $G_n$ tends to $1$ as $n \to \infty$.  
\end{corol}

We also obtain the following result on intersections of subgroups of almost simple exceptional groups as an immediate corollary of Theorem \ref{t:main}.

\begin{corol}\label{c:main}
Let $G$ be a finite almost simple exceptional group of Lie type with socle $G_0$ and let $H$ be any subgroup of $G$ not containing $G_0$. Then one of the following holds:
\begin{itemize}\addtolength{\itemsep}{0.2\baselineskip}
\item[{\rm (i)}] $G_0 = E_7(q)$ and $H$ is a subgroup of $P_7$. 
\item[{\rm (ii)}] $G_0 = E_6(q)$ and $H$ is a subgroup of $P_1$ or $P_6$.
\item[{\rm (iii)}] There exist $x_i \in G$, $i = 1, \ldots, 4$, such that 
$$H \cap H^{x_1} \cap H^{x_2} \cap H^{x_3} \cap H^{x_4} = 1.$$
\end{itemize}
\end{corol}

This is related to Problem 17.41(b) in the Kourovka notebook \cite{Kou}: \emph{Let $H$ be a solvable subgroup of a finite group $G$ that has no nontrivial solvable normal subgroups. Is it always possible to find $5$ conjugates of $H$ whose intersection is trivial?}

The example $G = S_8$, $H = S_4 \wr S_2$ shows that at least $5$ conjugates are needed, in general. Through work of Vdovin \cite{Vdo}, this problem has essentially been reduced to almost simple groups. Furthermore, for an almost simple exceptional group $G$, he has shown that the conclusion holds with $4$ conjugates when $H$ is a solvable Hall subgroup of $G$ (see \cite{Vdovin}). In some sense, part (iii) of Corollary \ref{c:main} settles the original problem (with $5$ conjugates) for almost \emph{all} subgroups of exceptional groups (if $G_0 = E_7(q)$, for example, the only possible exceptions are subgroups for which every maximal overgroup is a $P_7$ parabolic). See recent work of Baikalov \cite{Bai} for analogous results on solvable subgroups of symmetric and alternating groups.

\vs

Let us say a few words on the proof of Theorem \ref{t:main}. Let $G \leqs {\rm Sym}(\O)$ be a finite almost simple primitive non-standard permutation group with socle $G_0$ and point stabilizer $H$. Assume that $b(G)=6$. By the main theorems of \cite{Bur7,BGS} we immediately deduce that $G_0$ is either a sporadic group or an exceptional group of Lie type. In the former case, the examples arising in part (i) of Theorem \ref{t:main} are determined in \cite{BOB}, so we may assume that $G_0$ is an exceptional group. Here the main theorem of \cite{BLS} gives $b(G) \leqs 6$, so further work is needed to determine when equality holds. 

To proceed, we refine the probabilistic approach from \cite{BLS}, which relies on fixed point ratio estimates (this method for studying base sizes was originally introduced by Liebeck and Shalev in \cite{LSh2}). To describe the general set-up, fix a positive integer $c$ and let $Q(G,c)$ be the probability that a randomly chosen $c$-tuple of points in $\O$ is \emph{not} a base for $G$, so $b(G) \leqs c$ if and only if $Q(G,c)<1$. As explained in \cite{BLS}, it is easy to show that 
\begin{equation}\label{e:bd}
Q(G,c) \leqs \sum_{i=1}^{k}{|x_{i}^{G}|\cdot {\rm fpr}(x_{i})^{c}}=:\what{Q}(G,c),
\end{equation}
where $x_1, \ldots, x_k$ represent the distinct conjugacy classes in $G$ of elements of prime order. Here 
\begin{equation}\label{e:fpr}
{\rm fpr}(x) = \frac{|x^{G}\cap H|}{|x^{G}|}
\end{equation}
is the \emph{fixed point ratio} of $x \in G$, which is the proportion of points in $\O$ that are fixed by $x$. In view of \eqref{e:bd}, we can apply upper bounds on fixed point ratios to bound $\what{Q}(G,c)$ from above. 

Detailed results on fixed point ratios for exceptional groups of Lie type are presented in \cite{LLS2}, some of which rely on bounds on the dimensions of fixed point spaces for actions of the corresponding exceptional algebraic groups, which are determined in the companion paper \cite{LLS}. We will make extensive use of these results, together with improved estimates in several cases of interest. Our aim in almost every case is to identify a function $f(q)$, where $q$ is the size of the underlying field, such that 
$$\what{Q}(G,5) \leqs f(q) < 1$$
for all $q$, with the additional property that $f(q) \to 0$ as $q \to \infty$ (the latter property is needed for Corollary \ref{t:main2}). 

As in \cite{BLS}, we make a distinction between the cases where $H$ is a parabolic or non-parabolic subgroup of $G$. Following \cite{BLS} and \cite{LLS2}, we will use a character-theoretic approach to handle the parabolic actions, which relies in part on L\"{u}beck's work \cite{Lub} on Green functions for exceptional groups of Lie type; see Section \ref{s:para} for further details. For non-parabolic actions, we first reduce the problem to a handful of cases of the form $H = N_G(\bar{H}_{\s})$, where $\bar{H}$ is a $\s$-stable positive dimensional closed subgroup of the ambient simple algebraic group $\bar{G}$ and $\s$ is an appropriate Steinberg endomorphism of $\bar{G}$. By studying the embedding of $\bar{H}$ in $\bar{G}$, we can derive bounds on $|x^G \cap H|$ for elements $x \in G$ of prime order, which yield upper bounds on ${\rm fpr}(x)$ via \eqref{e:fpr}. In this regard, work of Lawther \cite{Lawunip} is particularly useful when $x$ is a unipotent element.

\begin{remk}
It is natural to consider the possibility of extending Theorem \ref{t:main} to the non-standard groups $G$ with $b(G)=5$. Again, there is an immediate reduction to exceptional groups. Indeed, the main theorem of \cite{Bur7} implies that  
$$(G,H) = ({\rm U}_{6}(2).2, {\rm U}_{4}(3).2^2)$$ is the only non-standard classical group with $b(G)=5$. Similarly, if $G_0$ is an alternating group then 
$(G,H) = (S_6, {\rm PGL}_{2}(5))$
is the only example (see \cite{BGS}). If $G_0$ is sporadic then the main theorem of \cite{BOB} implies that $b(G)=5$ if and only if $(G,H)$ is one of the following:
\setlength\extrarowheight{3pt}
$$\begin{array}{llll}
({\rm M}_{12},{\rm M}_{11}) & ({\rm M}_{22}, {\rm L}_{3}(4)) & ({\rm M}_{22}.2, {\rm L}_{3}(4).2) & ({\rm HS}, {\rm M}_{22}) \\
({\rm HS}.2, {\rm M}_{22}.2) & ({\rm McL}, {\rm U}_{4}(3)) & ({\rm McL}.2, {\rm U}_{4}(3).2) & ({\rm Fi}_{22}, 2.{\rm U}_{6}(2)) \\
({\rm Fi}_{23}, 2.{\rm Fi}_{22}) & ({\rm Co}_{1}, {\rm Co}_{2}) & ({\rm Fi}_{24}', {\rm Fi}_{23}) & ({\rm Fi}_{24}, {\rm Fi}_{23} \times 2)
\end{array}$$
To handle the exceptional groups, it may be feasible to proceed as in the proof of Theorem \ref{t:main}. However, the analysis will be significantly more complicated (with many more cases to consider) and some substantial improvements to existing bounds on fixed point ratios will be required.
\end{remk}

\vs

Finally, let us comment on the notation used in this paper. Most of our group-theoretic notation is standard. In particular, we use the notation from \cite{KL} for simple groups of Lie type, so we write 
${\rm L}_{n}^{+}(q) = {\rm L}_{n}(q) = {\rm PSL}_{n}(q)$ and $E_6^{-}(q) = {}^2E_6(q)$, etc. We also use ${\rm P\O}_{n}^{\e}(q)$ to denote a simple orthogonal group, which differs from the notation in the Atlas \cite{Atlas}. In addition, $(a,b)$ denotes the greatest common divisor of integers $a$ and $b$, and we use $i_r(X)$ for the number of elements of order $r$ in a subset $X$ of a finite group. Further notation will be introduced as and when needed.

\section{Parabolic actions}\label{s:para}

Let $G \leqs {\rm Sym}(\O)$ be a finite almost simple primitive non-standard permutation group with socle $G_0$ and point stabilizer $H$. As explained in the introduction, in order to prove Theorem \ref{t:main} (as well as Corollaries \ref{t:main2} and \ref{c:main}) we may assume that $G$ is an exceptional group of Lie type over $\mathbb{F}_{q}$, where $q=p^f$ for a prime $p$. In this section, we assume $H$ is a maximal parabolic subgroup of $G$. 

Let $\bar{G}$ be the ambient simple algebraic group over the algebraic closure of $\mathbb{F}_{q}$ and let $\s$ be a Steinberg endomorphism of $\bar{G}$ such that $G_0 = (\bar{G}_{\s})'$. In addition, let $\bar{H}$ be a $\s$-stable parabolic subgroup of $\bar{G}$ corresponding to $H$. We begin with a reduction theorem (in the statement, $G$ is an almost simple exceptional group of Lie type with socle $G_0$).

\begin{thm}\label{t:par}
Suppose $H$ is a maximal parabolic subgroup of $G$. Then $b(G) \leqs 6$, with equality only if one of the following holds:
\begin{itemize}\addtolength{\itemsep}{0.2\baselineskip}
\item[{\rm (i)}] $G_0 = E_7(q)$ and $H=P_7$. 
\item[{\rm (ii)}] $G_0 = E_6(q)$ and $H=P_1$ or $P_6$.
\item[{\rm (iii)}] $G_0 = E_6(q)$, $G$ contains graph or graph-field automorphisms, and $H = P_{1,6}$.
\end{itemize}
\end{thm}

\begin{proof}
This follows immediately from \cite[Theorem 3]{BLS}. 
\end{proof}

In part (iii) we write $P_{1,6}$ to denote the intersection of appropriate parabolic subgroups of type $P_{1}$ and $P_{6}$ (here the Levi factor is of type $D_4$). Note that $P_{1,6}$ is maximal in $G$ if and only if $G$ contains graph or graph-field automorphisms of $G_0$.

\begin{prop}\label{p:par1}
We have $b(G)=6$ in cases (i) and (ii) of Theorem \ref{t:par}.
\end{prop}

\begin{proof}
Consider the natural action of $\bar{G}$ on the coset variety $\bar{G}/\bar{H}$ and let $b^0(\bar{G},\bar{H})$ be the connected base size of $\bar{G}$ as defined in \cite{BGS4} (this is the smallest integer $c$ such that $\bar{G}/\bar{H}$ contains $c$ points with finite pointwise stabilizer). By \cite[Theorem 6]{BGS4} we have $b^0(\bar{G},\bar{H})=6$, so 
\cite[Proposition 2.7(ii)]{BGS4} implies that $b(G) \geqs 6$ if $q>2$. In particular, we conclude that $b(G)=6$ if $q>2$. 

Finally, suppose $q=2$. In case (i), we observe that there is a non-empty open subset $U$ of $(\bar{G}/\bar{H})^5$ such that the stabilizer in $\bar{G}$ of every $5$-tuple in $U$ is $8$-dimensional. As explained in \cite[Remark 5.6]{BGS4}, this implies that $b(G)=6$. In (ii), a {\sc Magma} \cite{Magma} calculation shows that $b(G)=6$ (see \cite[Remark 1]{BLS}).
\end{proof}

In order to complete the proof of Theorem \ref{t:main} for parabolic actions, it remains to show that $b(G) \leqs 5$ in case (iii) of Theorem \ref{t:par}.

\begin{prop}\label{p:e6p16}
We have $b(G) \leqs 5$ in case (iii) of Theorem \ref{t:par}.
\end{prop}

\begin{proof}
Here $G_0=E_6(q)$ and $H = P_{1,6}$. Define $\what{Q}(G,5)$ as in \eqref{e:bd} and write 
\begin{equation}\label{e:alpha}
\what{Q}(G,5) = \a+\b+\gamma,
\end{equation}
where $\alpha$ and $\beta$ are the contributions to $\what{Q}(G,5)$ from semisimple and unipotent elements in $G$, respectively, and $\gamma$ is the contribution from graph, field and graph-field automorphisms of $G_0$. To estimate $\a$ and $\b$ we work with the corresponding permutation character $\chi = 1_{P_{1,6}}^{\bar{G}_{\s}}$, so that ${\rm fpr}(x) = \chi(x)/|\O|$ for all $x \in \bar{G}_{\s}$. Let $W \cong {\rm U}_{4}(2).2$ be the Weyl group of $\bar{G}$ and note that 
\begin{equation}\label{e:om}
|\O| = \frac{(q^9-1)(q^8+q^4+1)(q^5-1)(q^4+1)}{(q-1)^2}.
\end{equation}

We start by estimating $\a$. To do this, we proceed as in \cite[Section 3.1(ii)]{BLS}, using \cite[Corollary 3.2]{LLS2} to compute $\chi(x)$ for any semisimple element $x \in G$. In order to explain the set-up, fix a set $\Pi = \{\a_1, \ldots, \a_6\}$ of simple roots for $\G$ and let $\a_0$ be the highest root in the root system of $\bar{G}$. Recall that the semisimple classes in $\bar{G}_{\s}$ are parametrized by pairs $(J,[w])$, where $J$ is a proper subset of $\Pi \cup \{\a_0\}$ (determined up to $W$-conjugacy), $W_J$ is the subgroup of $W$ generated by the reflections in the roots in $J$, and $[w] = W_Jw$ is a conjugacy class representative of $N_W(W_J)/W_J$. For $\bar{G}_{\s} = E_6(q)$, the $\bar{G}_{\s}$-classes of $\s$-stable maximal tori of $\bar{G}$ are parametrized by the conjugacy classes $C_1, \ldots, C_k$ of $W$. Let $T_i$ be a maximal torus of $\bar{G}_{\s}$ corresponding to $C_i$ and let $r_i$ be the relative rank of $T_i$ (that is, $r_i$ is the multiplicity of $q-1$ as a divisor of $|T_i|$). Finally, let $W_{\bar{H}} \leqs W$ be the Weyl group of the corresponding parabolic subgroup $\bar{H}$ of $\bar{G}$. 

If $x \in G$ corresponds to the pair $(J,[w])$ as above, then \cite[Corollary 3.2]{LLS2} gives
\begin{equation}\label{e:c32}
\chi(x) = \sum_{i=1}^{k}\frac{|W|}{|C_i|}\cdot \frac{|W_{\bar{H}} \cap C_i|}{|W_{\bar{H}}|}\cdot \frac{|W_Jw \cap C_i|}{|W_J|}\cdot \frac{(-1)^{r_i}|(C_{\bar{G}}(x)^{0})_{\s}|_{p'}}{|T_i|},
\end{equation}
which can be evaluated using {\sc Magma} \cite{Magma}.  For example, suppose $q$ is odd and $x \in G$ is an involution with $C_{\bar{G}}(x) = A_5A_1$. Here $J = \{\a_0, \a_1, \a_3, \a_4, \a_5, \a_6\}$, $w=1$ and we compute
\begin{align*}
\chi(x) = & \, q^{12} + 2q^{11} + 8q^{10} + 16q^{9} + 26q^{8} + 33q^{7} + 41q^{6} + 42q^{5} + 38q^{4} + 28q^3 \\
& + 20q^2 + 11q + 4,
\end{align*}
which gives ${\rm fpr}(x)<q^{-11}$. Since there are fewer than $2q^{40}$ conjugates of $x$ in $G$, we deduce that the contribution to $\a$ from these elements is less than $2q^{40}(q^{-11})^5 = 2q^{-15}$. In this way, by repeatedly applying \eqref{e:c32}, we calculate that $\a<q^{-5}$.

Next let us turn to $\b$. Here we adopt the approach of \cite[Section 3.1(i)]{BLS}, which is based on the method originally introduced in \cite[Section 2]{LLS2}. As recorded in \cite[p.123]{BLS}, we can write 
$$\chi = R_{\phi_{1,0}}+2R_{\phi_{6,1}}+3R_{\phi_{20,2}}+R_{\phi_{15,5}}+R_{\phi_{30,3}}+2R_{\phi_{64,4}}+R_{\phi_{24,6}}$$
where each $R_{\phi}$ is an almost character of $\bar{G}_{\s}$ labelled by a complex irreducible character $\phi$ of $W$ (we use the labelling given in \cite[Section 13.2]{Carter}). The restriction of the $R_{\phi}$ to unipotent elements gives the Green functions of $\bar{G}_{\s}$, which are polynomials in $q$ with non-negative coefficients. By implementing an algorithm of Lusztig \cite{Lus}, the relevant Green functions have been computed by L\"{u}beck \cite{Lub}, which allows us to calculate $\chi(x)$ for each element $x \in \bar{G}_{\s}$ of order $p$. (In the general setting, there is a problem with certain unspecified scalars that arises in the implementation of Lusztig's algorithm, but this issue does not occur when $\bar{G}=E_6$ and $x$ has order $p$; see \cite[Section 3.1(i)]{BLS} for further discussion.) 
In addition, we can read off $|x^{G}|$ from \cite[Table 22.2.3]{LSbook} and this allows us to compute $\b$ precisely.

For instance, if $x \in G$ is a long root element then 
\begin{align*}
\chi(x) = & \, q^{18} + 2q^{17} + 5q^{16} + 7q^{15} + 11q^{14} + 14q^{13} + 18q^{12} + 18q^{11} + 19q^{10} + 18q^9 \\
& + 17q^8 + 13q^7 + 11q^6 + 9q^5 + 7q^4 + 4q^3 + 3q^2 + 2q + 1 \\
|x^G| = & \, (q^3+1)(q^4+1)(q^6+1)(q^9-1)
\end{align*}
and we deduce that the contribution to $\b$ from these elements is less than $q^{-5}$. In fact, by computing $\chi(x)$ for all $x \in \bar{G}_{\s}$ of order $p$, we see that 
$\b<q^{-5}$.

To complete the proof, it remains to consider the contribution $\gamma$ from graph, field and graph-field automorphisms of $G_0$. First assume that $x \in G$ is an involutory graph automorphism. If $C_{\bar{G}}(x) = F_4$ then ${\rm fpr}(x) \leqs q^{-9}$ if $q$ is even (see \cite[Section 3.1(iv)]{BLS}) and by arguing as in the proof of \cite[Lemma 6.4]{LLS2} we deduce that the same bound also holds when $q$ is odd.  Similarly, if $C_{\bar{G}}(x) \ne F_4$ then ${\rm fpr}(x) \leqs q^{-13}$ for all $q$. It follows that the contribution from graph automorphisms is less than
$$2q^{26}(q^{-9})^5+2q^{42}(q^{-13})^5<q^{-17}.$$
Similarly, if $x$ is an involutory graph-field automorphism then $q=q_0^2$ and 
$|C_{\O}(x)|$ is equal to the index of a $P_{1,6}$ parabolic subgroup in ${}^2E_6(q_0)$ (see \cite[p.452]{LLS2}). This implies that ${\rm fpr}(x)<q^{-12}$, so the contribution from these elements is at most $6q^{39}(q^{-12})^5 <q^{-19}$.

Finally, suppose $x \in G$ is a field automorphism of prime order $r$. Here $q=q_0^r$ and we observe that
$$|x^G|<2q^{78(1-r^{-1})} = g(r,q),\;\; {\rm fpr}(x) = f(q_0)/f(q) = h(r,q),$$
where $f(q) = |\O|$ as in \eqref{e:om}. Set $j(r,q) = g(r,q)h(r,q)^5$ and let $\pi$ be the set of prime divisors of $\log_{p}q$. Then the contribution to $\what{Q}(G,5)$ from field automorphisms is less than 
$$\sum_{r \in \pi}(r-1)\cdot j(r,q)<j(2,q)+2j(3,q)+4j(5,q)+\log_{2}q.q^{78}h(7,q)^5<q^{-12}$$
and thus $\gamma<q^{-12}+q^{-17}+q^{-19} < q^{-11}$.

Therefore, by bringing together the above bounds, we conclude that 
$$\what{Q}(G,5)< 2q^{-5}+q^{-11}$$
and the desired result follows.
\end{proof}

\begin{rem}\label{r:par}
Notice that in the proof of the previous proposition, we identify an explicit function $f(q)$ such that $\what{Q}(G,5)<f(q)$ and $f(q) \to 0$ as $q \to \infty$. Therefore, in this particular case, the probability that $5$ randomly chosen points in $\O$ form a base for $G$ tends to $1$ as $|G|$ tends to infinity. By combining this observation with \cite[Theorem 3]{BLS}, we see that the same conclusion holds whenever $G_0$ is an exceptional group and $H$ is a maximal parabolic subgroup, with the possible exception of the cases recorded in parts (i) and (ii) of Theorem \ref{t:par}. Indeed, these cases are genuine exceptions by Proposition \ref{p:par1}, which explains the additional condition on the sequence $(G_n)$ in Corollary \ref{t:main2}.
\end{rem}

\section{Non-parabolic actions}

For the remainder of the paper, we may assume that $G$ is an exceptional group of Lie type and the point stabilizer $H$ is a maximal non-parabolic subgroup. As in the previous section,  our aim is to establish the bound $\what{Q}(G,5)<1$ through a careful analysis of fixed point ratios for prime order elements. In addition, we will show that $\what{Q}(G,5) \to 0$ as $q \to \infty$, which will allow us to complete the proofs of Theorem \ref{t:main} and Corollary \ref{t:main2}. In view of \cite[Theorem 4]{BLS}, we may assume that 
$$G_0 \in \{E_7(q), E_{6}^{\e}(q), F_4(q)\}.$$

Let us fix some additional notation. As before, let $\bar{G}$ be the ambient simple algebraic group (of type $E_7$, $E_6$ or $F_4$) over the algebraic closure of $\mathbb{F}_{q}$, so $G_0 = (\bar{G}_{\s})'$ for an appropriate Steinberg endomorphism $\s$. For $x \in \bar{G}$ and a closed subgroup $\bar{H}$ we set
\begin{equation}\label{e:dbar}
\bar{D} = C_{\bar{G}}(x),\;\; \delta(x) = \dim x^{\bar{G}} - \dim (x^{\bar{G}} \cap \bar{H}).
\end{equation}
We will repeatedly use standard information on semisimple and unipotent conjugacy classes in $\bar{G}_{\s}$ (including class lengths, etc.) recorded in \cite{FJ} (semisimple classes in $E_6^{\e}(q)$ and $E_7(q)$), \cite{Shinoda, Shoji} (conjugacy classes in $F_4(q)$) and \cite{LSbook} (unipotent classes in all cases). In particular, we will adopt the notation from \cite{LSbook} for unipotent classes in $\bar{G}$. We also refer the reader to L\"{u}beck's online data \cite{Lu2}, which provides detailed information on semisimple conjugacy classes in exceptional groups.

As in the proof of Proposition \ref{p:e6p16}, it will be convenient to write 
$$\what{Q}(G,5) = \a+\b+\gamma,$$ 
where $\a$ is the contribution from semisimple elements, $\b$ from unipotent elements and $\gamma$ from field, graph and graph-field automorphisms. In order to estimate these quantities, we will frequently apply the following elementary result (see \cite[Proposition 2.3]{BLS}).

\begin{prop}\label{p:bdd}
Let $G \leqs {\rm Sym}(\O)$ be a finite transitive permutation group with point stabilizer $H$. Suppose $x_1, \ldots, x_m$ represent distinct $G$-classes such that $\sum_{i}|x_i^G \cap H| \leqs A$ and $|x_i^G|\geqs B$ for all $i$. Then
$$\sum_{i=1}^{m}|x_i^G|\cdot {\rm fpr}(x_i)^c \leqs B(A/B)^c$$
for any positive integer $c$.
\end{prop}

The next four sections are dedicated to the proof of Theorem \ref{t:main} in the cases where $G_0 = E_7(q)$, $E_6^{\e}(q)$ with $q \geqs 3$, $E_6^{\e}(2)$ and $F_4(q)$, respectively. We complete the proof of Corollary \ref{t:main2} in Section \ref{ss:c2}.

\subsection{$G_0 = E_7(q)$}\label{s:e7}

In this section we handle the case $G_0 = E_7(q)$. For $q=2$, the possibilities for $H$ are determined in \cite{BBR}. For information on semisimple and unipotent conjugacy classes, we refer the reader to \cite[Section 3]{FJ} and \cite[Table 22.2.2]{LSbook}, respectively. We begin by reducing the problem to subgroups $H$ with $|H|>q^{46}$.

\begin{lem}\label{l:e7_bd}
If $G_0 = E_7(q)$ and $|H| \leqs q^{46}$, then $b(G) \leqs 5$.
\end{lem}

\begin{proof}
As explained in the proof of \cite[Lemma 4.8]{BLS}, we have
\begin{equation}\label{e:1}
\what{Q}(G,5) < b(a/b)^5 + cd^5
\end{equation}
where
$a = q^{46}$, $b=\frac{1}{2}(q+1)^{-1}q^{65}$, $c=3q^{55}$ and $d=2q^{-12}$. Therefore, if $q \geqs 3$ we have $\what{Q}(G,5)<2q^{-1}$ and the result follows. For the remainder, let us assume $q=2$, so $G = E_7(2)$ and the possibilities for $H$ are given in \cite[Theorem 1.1]{BBR}. Let $x \in G$ be an element of prime order $r$. 

We claim that ${\rm fpr}(x) \leqs 2^{-12}$ if $\dim x^{\bar{G}}<64$, which implies that \eqref{e:1} holds with $d=2^{-12}$ (and $a,b,c$ as above) and thus $\what{Q}(G,5)<1$ as required. In view of \cite[Theorem 2]{LLS2}, we may assume that $r=2$ and $H = N_G(\bar{H}_{\s})$ is a subgroup of maximal rank. Since $|H| \leqs q^{46}$, it follows that 
$$\mbox{$\bar{H}^0=A_2A_5$, $A_1^3D_4$, $A_1^7$ or $T_7$.}$$
(Here, and elsewhere, we write $T_i$ for an $i$-dimensional torus.) In addition, the condition $\dim x^{\bar{G}}<64$ implies that $x$ belongs to one of the classes labelled $A_1$, $A_1^2$ or $(A_1^3)^{(1)}$ in \cite[Table 22.2.2]{LSbook}.  Note that $|x^G|>2^{34}$.

First assume $\bar{H}^0 = T_7$, in which case $H = 3^7.(2 \times {\rm Sp}_{6}(2))$ and $|H|<2^{33}$. If $x$ is not in the class $A_1$, then $|x^G|>2^{52}$ and the desired bound holds. On the other hand, if $x$ is a long root element then by arguing as in the proof of \cite[Lemma 4.3]{LLS2}, noting that the Weyl group $W(E_7)$ contains $63$ reflections, we deduce that $|x^G \cap H| \leqs 3^7.63 < 2^{18}$ and the result follows from the bound $|x^G|>2^{34}$. Similarly, if $\bar{H}^0 = A_1^7$ then $H = {\rm L}_{2}(2^7).7$ and $i_2(H)=2^{14}-1$, so 
$${\rm fpr}(x) = \frac{|x^G \cap H|}{|x^G|} < \frac{2^{14}-1}{2^{34}} < 2^{-20}$$
and the claim holds. 

Next assume $\bar{H}^0=A_2A_5$. Here $H = ({\rm SL}_{3}^{\e}(2) \circ {\rm SL}_{6}^{\e}(2)).A$, with $A = Z_2$ ($\e=+$) or $S_3$ ($\e=-$). If $x$ is in $A_1$ then by arguing as in the proof of \cite[Lemma 4.4]{LLS2}, using  \cite[Proposition 1.13(ii)]{LLS}, we calculate that 
$${\rm fpr}(x) \leqs \frac{|x_1^{{\rm SL}_{3}^{\e}(2)}|+|x_2^{{\rm SL}_{6}^{\e}(2)}|}{|x^G|}<2^{-23},$$
where $x_1$ and $x_2$ are long root elements in ${\rm SL}_{3}^{\e}(2)$ and ${\rm SL}_{6}^{\e}(2)$, respectively, and $\e=\pm$ according to the specific structure of $H$. Similarly, if $x$ is in  $A_1^2$ or $(A_1^3)^{(1)}$ then $|x^G|>2^{52}$ and thus ${\rm fpr}(x)<2^{-25}$ since 
$$i_2(H) \leqs i_2({\rm SL}_{3}^{\e}(2).2 \times {\rm SL}_{6}^{\e}(2).2)<2^{27}.$$

Finally, suppose that $\bar{H}^0=A_1^3D_4$, so 
$H =({\rm L}_{2}(8) \times {}^3D_4(2)).3$ or $((S_3)^3 \times \O_{8}^{+}(2)).S_3$.
If $x$ is in $A_1^2$ or $(A_1^3)^{(1)}$ then ${\rm fpr}(x)<2^{-52}|H|<2^{-12}$, so we may assume $x$ is a long root element. If $H = ({\rm L}_{2}(8) \times {}^3D_4(2)).3$ then \cite[Proposition 1.13(ii)]{LLS} implies that 
$${\rm fpr}(x) \leqs \frac{|x_1^{{}^3D_4(2)}|}{|x^G|}< \frac{819}{2^{34}}<2^{-24},$$
where $x_1 \in {}^3D_4(2)$ is a long root element. The case $H = ((S_3)^3 \times \O_{8}^{+}(2)).S_3$ is very similar.

This justifies the claim and we conclude that $b(G) \leqs 5$, as required.
\end{proof}

To complete the proof of Theorem \ref{t:main} for $G_0 = E_7(q)$, we may assume that $|H|>q^{46}$, in which case \cite[Lemma 4.7]{BLS} implies that either
\begin{itemize}\addtolength{\itemsep}{0.2\baselineskip}
\item[{\rm (i)}] $H = N_G(\bar{H}_{\s})$ with $\bar{H}^0 = A_1F_4$, $A_1D_6$, $A_7$ or $T_1E_6$; or
\item[{\rm (ii)}] $q=q_0^2$ and $H$ is a subfield subgroup of type $E_7(q_0)$.
\end{itemize}
In (ii), the proof of \cite[Proposition 4.11]{BLS} gives $b(G) \leqs 4$. Similarly, $b(G) \leqs 5$ if $\bar{H}^0 = A_1F_4$ in (i) (see the proof of \cite[Lemma 4.10]{BLS}). For the three remaining possibilities, \cite{BLS} only gives $b(G) \leqs 6$, so we need to improve the bound in each of these cases.

It will be convenient to introduce some additional terminology. For a semisimple element $x \in G$, let us say that $\bar{D} = C_{\bar{G}}(x)$ is \emph{large} if $\bar{D}^0$ has an $E_6$, $D_6$ or $A_7$ factor (otherwise $\bar{D}$ is \emph{small}), and write $\a=\a_1+\a_2$, where $\a_1$ is the contribution to $\what{Q}(G,5)$ from the semisimple elements with large centralizer. 

\begin{lem}\label{l:e7a1d6}
If $G_0 = E_7(q)$ and $H = N_G(\bar{H}_{\sigma})$ with $\bar{H} = A_1D_6$, then $b(G) \leqs 5$.
\end{lem}

\begin{proof}
Write $\what{Q}(G,5) = \a+\b+\gamma$ as in \eqref{e:alpha}. By arguing as in the proof of \cite[Lemma 4.10]{BLS}, we have
$$\gamma<\log_2q.q^{133}(8q^{-32})^5<q^{-11}.$$

Next we claim that $\a<q^{-16}$. If $q \geqs 3$ then the proof of \cite[Lemma 4.10]{BLS} yields 
$$\a<q^{133}(q^{-30})^5+q^{71}(q^{-19})^5<q^{-16}$$
as required. Now assume $q=2$, so $G=E_7(2)$ and $H = S_3 \times \O_{12}^{+}(2)$. Note that the semisimple classes in $G$ are listed in \cite[Table 2]{BBR}, together with the structure of the corresponding centralizers. 

Let $x \in G$ be an element of odd prime order $r$ and define $\bar{D}$ and $\delta(x)$ as in \eqref{e:dbar}.  First assume $\bar{D}$ is large, in which case $r=3$ (more precisely, $x$ is a $3A$ or $3B$ element in the notation of \cite[Table 2]{BBR}), $|x^G|>\frac{1}{2}2^{\dim x^{\bar{G}}}$ and \cite[Theorem 2]{LLS} gives $\delta(x) \geqs 24$. Now $H$ has $15$ classes of elements of order $3$ and one checks that $|y^H|<4 \cdot 2^{\dim y^{\bar{H}}}$ for all $y \in H$ of order $3$. Therefore, by arguing as in the proof of \cite[Lemma 4.5]{LLS2}, we deduce that
$${\rm fpr}(x) <\frac{15 \cdot 4 \cdot 2^{\dim(x^{\bar{G}} \cap \bar{H})}}{\frac{1}{2}2^{\dim x^{\bar{G}}}} = 120\cdot 2^{-\delta(x)} <2^{-17}.$$
Since $G$ contains at most $2^{66}$ such elements, it follows that $\a_1< 2^{66}(2^{-17})^5 = 2^{-19}$.

Now assume $\bar{D}$ is small. Here $r \in \{3,5,7,11,17,31\}$, $\delta(x) \geqs 40$ and by arguing as above we deduce that  
${\rm fpr}(x) < 120\cdot 2^{-\delta(x)} <2^{-33}$.
Therefore, $\a_2<2^{133}(2^{-33})^5 = 2^{-32}$ and we conclude that $\a<q^{-16}$ for all $q \geqs 2$.

Finally, let $x \in G$ be a unipotent element of order $p$. The $\bar{G}$-class of each unipotent element in $\bar{H}$ is determined in \cite[Section 4.10]{Lawunip} and this allows us to compute $\delta(x)$ precisely, which yields a good upper bound on ${\rm fpr}(x)$. First assume $p$ is odd. As explained in the proof of \cite[Lemma 4.10]{BLS}, if $\dim x^{\G}>66$ then ${\rm fpr}(x)<q^{-26}$, which implies that the contribution to $\b$ from these elements is less than $q^{126}(q^{-26})^5=q^{-4}$ (here we are using the fact that $G$ has precisely $q^{126}$ unipotent elements). Similarly, if $\dim x^{\G} \leqs 66$ then one checks that ${\rm fpr}(x)<2q^{-16}$. For example, if $x$ is a long root element then $|x^G|>q^{34}$ and
$$|x^G \cap H| \leqs |x_1^{{\rm SL}_{2}(q)}| + |x_2^{\O_{12}^{+}(q)}|<2q^{18},$$
where $x_1$ and $x_2$ are long root elements in ${\rm SL}_{2}(q)$ and $\O_{12}^{+}(q)$, respectively, whence ${\rm fpr}(x)<2q^{-16}$ as claimed. Since $G$ contains at most $2q^{66}$ such elements, it follows that 
$$\b<q^{-4}+2q^{66}(2q^{-16})^5<q^{-3}.$$ 
A similar calculation when $p=2$ also shows that $\b<q^{-3}$ (here the $\bar{G}$-class of each involution in $\bar{H}$ is recorded in \cite[Table 3]{BLS}).

We conclude that
$$\what{Q}(G,5)< q^{-3}+q^{-11}+q^{-16}$$
and the result follows.
\end{proof}

\begin{lem}\label{l:e7a7}
If $G_0 = E_7(q)$ and $H = N_G(\bar{H}_{\sigma})$ with $\bar{H} = A_7.2$, then $b(G) \leqs 5$.
\end{lem}

\begin{proof}
Define $\a,\b,\gamma,\a_1$ and $\a_2$ as before. We start by estimating $\a$. 

First assume $q \geqs 3$. Let $x \in G$ be a semisimple element of prime order $r$. If $\bar{D}$ is large then ${\rm fpr}(x) \leqs q^{-19}$ by \cite[Theorem 2]{LLS2} and we calculate that $G$ has at most $q^{71}$ such elements. On the other hand, if $\bar{D}$ is small then by appealing to \cite[Lemma 4.5]{LLS2} we deduce that
$${\rm fpr}(x) <  8|W(E_7):W(A_7)|\left(\frac{q+1}{q-1}\right)^7q^{-\delta(x)} \leqs 576\left(\frac{q+1}{q-1}\right)^7q^{-44}$$
since $\delta(x) \geqs 44$ by \cite[Theorem 2]{LLS}. Therefore
$$\a<q^{71}(q^{-19})^5+q^{133}\left(576\left(\frac{q+1}{q-1}\right)^7q^{-44}\right)^5<q^{-23}.$$

Now assume $q=2$, so $G=E_7(2)$ and $H = {\rm L}_{8}^{\e}(2).2$ with $\e=\pm$. Suppose $\bar{D}$ is large. As noted in the proof of the previous lemma, the condition on $\bar{D}$ implies that $r=3$ and there are fewer than $2^{66}$ such elements in $G$. Now ${\rm L}_{8}^{\e}(2)$ has $9-5\e$ classes of elements of order $3$; if $y \in H$ is such an element then $|y^{{\rm L}_{8}^{\e}(2)}|<4\cdot 2^{\dim y^{\bar{H}}}$ and $|y^G|>\frac{1}{2}2^{\dim y^{\bar{G}}}$. Therefore
${\rm fpr}(x) <112\cdot 2^{-\delta(x)}$
with $\delta(x) \geqs 27$ (see \cite[Theorem 2]{LLS}), so ${\rm fpr}(x) < 2^{-20}$. It follows that $\a_1<2^{66}(2^{-20})^5 = 2^{-34}$. Similarly, if $\bar{D}$ is small then $\delta(x) \geqs 44$ and 
${\rm fpr}(x) <144\cdot 2^{-\delta(x)} < 2^{-36}$, whence $\a_2<2^{133}(2^{-36})^5 = 2^{-47}$. In particular, $\a<q^{-23}$ for all $q \geqs 2$.

Next we consider $\b$. Let $x \in G$ be a unipotent element of order $p$. First assume $p$ is odd, so $x \in \bar{H}^0$ and we can compute $\delta(x)$ using the fusion information in \cite[Section 4.11]{Lawunip}. We find that $x^{\bar{G}} \cap \bar{H}$ is a union of at most two $\bar{H}$-classes, so 
the proof of \cite[Lemma 4.5]{LLS2} yields
\begin{equation}\label{e:nee}
{\rm fpr}(x)<48\left(\frac{q+1}{q-1}\right)^7q^{-\delta(x)}
\end{equation}
since $|C_{\bar{G}}(x):C_{\bar{G}}(x)^0| \leqs 6$ (see \cite[Table 22.1.2]{LSbook}). 
If $\dim x^{\bar{G}}>66$ then $\delta(x) \geqs 42$ and thus the contribution to $\b$ from these elements is at most
$$q^{126}\left(48\left(\frac{q+1}{q-1}\right)^7q^{-42}\right)^5<q^{-44}.$$
If $\dim x^{\bar{G}} \leqs 66$ then $x$ belongs to one of the $\bar{G}$-classes labelled $A_1$, $A_1^2$, $A_2$, $(A_1^3)^{(1)}$ or $(A_1^3)^{(2)}$. As noted in the proof of the previous lemma, there are fewer than $2q^{66}$ such elements in $G$ and we calculate that ${\rm fpr}(x)<2q^{-20}$. For example, if $x \not\in A_1$ then $\delta(x) \geqs 28$ and the bound on ${\rm fpr}(x)$ follows from \eqref{e:nee}. On the other hand, if $x \in A_1$ then $|x^G|>q^{34}$ and 
$$|x^G \cap H|  \leqs |x_1^{{\rm SL}_{8}^{\e}(q)}| = \frac{|{\rm SL}_{8}^{\e}(q)|}{q^{13}|{\rm GL}_{6}^{\e}(q)|}<2q^{14},$$
where $x_1 \in {\rm SL}_{8}^{\e}(q)$ is a long root element, which gives ${\rm fpr}(x)<2q^{-20}$. Therefore, the contribution from these elements is less than $2q^{66}(2q^{-20})^5 = 2^6q^{-34}$ and we conclude that $\b<q^{-44} + 2^6q^{-34} < q^{-30}$ when $p$ is odd.

Now suppose $p=2$. There are six conjugacy classes of involutions in $\bar{H} = A_7.2$; four classes in the connected component $\bar{H}^0=A_7$, with representatives $x_i$, $1 \leqs i \leqs 4$ (where $x_i$ has Jordan form $[J_2^i,J_1^{8-2i}]$ on the natural $A_7$-module), plus two classes of involutory graph automorphisms in $\bar{H} \setminus \bar{H}^0$. The latter two classes are represented by elements $\tau_1$ and $\tau_2$, where $C_{A_7}(\tau_1) = C_4$ and  $C_{A_7}(\tau_2) = C_{C_4}(t)$ with $t \in C_4$ a long root element. By inspecting \cite[Section 4.11]{Lawunip} and the proof of \cite[Lemma 3.18]{BGS4}, we deduce that $x_1$ is in the $\bar{G}$-class $A_1$, $x_2$ is in $A_1^2$, $x_3$ and $x_4$ are in $(A_1^3)^{(2)}$, $\tau_1$ is in $(A_1^3)^{(1)}$ and $\tau_2$ is in $A_1^4$. It follows that ${\rm fpr}(x)<2q^{-20}$ and thus $\b<2q^{70}(2q^{-20})^5=2^6q^{-30}$ since $G_0$ contains at most $2q^{70}$ involutions. We conclude that $\b<q^{-24}$ for all $q \geqs 2$.

Finally, let us estimate $\gamma$, which is the contribution to $\what{Q}(G,5)$ from field automorphisms. Let $x \in G$ be a field automorphism of prime order $r$, so $q=q_0^r$. If $r=2$ then \cite[Theorem 2]{LLS2} yields ${\rm fpr}(x) \leqs q^{-22}$, and we note that $G$ has at most $4q^{133/2}$ such elements. Similarly, if $r$ is odd then
$${\rm fpr}(x) \leqs \frac{2|{\rm SL}_{8}^{\e}(q):{\rm SL}_{8}^{\e}(q_0)|}{|E_7(q):E_7(q_0)|}<8q^{-70(1-r^{-1})} \leqs 8q^{-140/3}$$
and thus
$$\gamma<4q^{133/2}(q^{-22})^5+\log_2q.q^{133}(8q^{-140/3})^5<q^{-41}.$$

By combining the above estimates, we conclude that $\what{Q}(G,5)<q^{-23}+q^{-24}+q^{-41}$ and the proof of the lemma is complete.
\end{proof}

\begin{lem}\label{l:e7t1e6}
If $G_0 = E_7(q)$ and $H= N_G(\bar{H}_{\sigma})$ with $\bar{H} = T_1E_6.2$, then $b(G) \leqs 5$.
\end{lem}

\begin{proof}
Write $\what{Q}(G,5)=\a+\b+\gamma$ and $\a=\a_1+\a_2$ as before. As explained in the proof of \cite[Lemma 4.9]{BLS}, we have 
$$\gamma<2q^{133/2}(q^{-22})^5 + \log_2q.q^{133}(16q^{-36})^5 < q^{-26}.$$

Next we claim that $\a<q^{-2}$. To see this, let $x \in G$ be a semisimple element of prime order $r$ and let us assume $q \geqs 3$ for now. As in the proofs of the previous two lemmas, we have $\a_1<q^{71}(q^{-19})^5 = q^{-24}$. Now assume $\bar{D}$ is small, so $r$ is odd and \cite[(4.4)]{BLS} implies that
\begin{equation}\label{e:ner}
{\rm fpr}(x)<\frac{224(q+1)^z}{q^{\delta(x)+z-6}(q-1)^6}
\end{equation}
where $z = \dim Z(\bar{D}^0)$ and $\delta(x) \geqs 34$ (see \cite[Theorem 2]{LLS}). If $\delta(x) = 34$ then $\dim x^{\bar{G}} \leqs 106$ and by inspecting \cite{FJ} we calculate that $G$ has at most $q^{110}$ such elements. Therefore, 
$$\a_2<q^{133}\left(\frac{224(q+1)^7}{q^{37}(q-1)^6}\right)^5+q^{110}\left(\frac{224(q+1)^7}{q^{35}(q-1)^6}\right)^5$$
and we deduce that $\a<q^{-2}$.

Now assume $q=2$, so $G=E_7(2)$ and $H = E_6(2).2$ or $3.{}^2E_6(2).S_3$. As noted in the proof of \cite[Lemma 4.9]{BLS}, if $\bar{D}$ is large then ${\rm fpr}(x)<2^{-20}$ and we deduce that $\a_1<2^{69}(2^{-20})^5 = 2^{-31}$. Now assume $\bar{D}$ is small. If $r=3$ then 
$\bar{D}^0 = A_2A_5$, $A_6T_1$ or $A_1D_5T_1$ (see \cite[Table 4.7.1]{GLS}), $z \leqs 1$ and $\delta(x) \geqs 34$ as before, so \eqref{e:ner} implies that ${\rm fpr}(x)<2^{-19}$. Since $G$ contains fewer than $2^{91}$ elements of order $3$ (see \cite[Table 2]{BBR}), their contribution is at most $2^{91}(2^{-19})^5 = 2^{-4}$. Finally, suppose $r \geqs 5$. By considering the conjugacy classes in $E_6^{\e}(2)$, we see that $x^G \cap H$ is a union of at most $8$ distinct $H$-classes, so 
$${\rm fpr}(x)<\frac{8 \cdot 2 \cdot 2^{\dim (x^{\bar{G}} \cap \bar{H})}}{\frac{1}{2}2^{\dim x^{\bar{G}}}} = 32\cdot 2^{-\delta(x)} \leqs 2^{-29}.$$
Therefore, $\a<2^{-31}+2^{-4}+2^{133}(2^{-29})^5<2^{-3}$.

Finally, let us turn to $\b$. If $p=2$ then 
$\b<\sum_{i=1}^{5}c_id_i^5<q^{-19}$, where the terms $c_i,d_i$ are presented in \cite[Table 8]{BLS}. Now assume $p$ is odd and let $x \in G$ be a unipotent element of order $p$. As explained in the proof of \cite[Lemma 4.9]{BLS}, we have 
\begin{equation}\label{e:nej}
{\rm fpr}(x)<\frac{24(q+1)^7}{q^{\delta(x)+1}(q-1)^6}
\end{equation}
and we can calculate $\delta(x)$ using \cite{Law} and the restriction $V_{56}{\downarrow}E_6 = V_{27} \oplus (V_{27})^* \oplus 0^2$ of the minimal module $V_{56}$ for $E_7$ (here $V_{27}$ is one of the irreducible $27$-dimensional modules for $E_6$, and $0$ denotes the trivial module). 

First assume $\dim x^{\bar{G}}>66$, in which case $\delta(x) \geqs 30$. If $p=3$ then \eqref{e:nej} implies that ${\rm fpr}(x)<q^{-23}$ and thus the contribution to $\b$ from these elements is less than $q^{92}(q^{-23})^5 = q^{-23}$ since $i_3(G_0)<q^{92}$ by \cite[Proposition 1.3(ii)]{LLS2}. Similarly, if $p \geqs 5$ then ${\rm fpr}(x)<q^{-26}$ and the contribution is less than $q^{126}(q^{-26})^5 = q^{-4}$. Next assume $34<\dim x^{\bar{G}} \leqs 66$. By considering each $\bar{G}$-class in turn, we calculate that ${\rm fpr}(x) < q^{-18}$, so the contribution is at most $2q^{66}(q^{-18})^5 = 2q^{-24}$. Finally, the contribution from long root elements is less than $2q^{34}(2q^{-12})^5 = 2^6q^{-26}$ since ${\rm fpr}(x) \leqs 2q^{-12}$ by \cite[Theorem 2]{LLS2}. We conclude that $\b<q^{-3}$ for all $q$. 

By bringing together the above bounds, we see that 
$\what{Q}(G,5)<q^{-2} + q^{-3} + q^{-26}$
and the result follows.
\end{proof}

\vs 

This completes the proof of Theorem \ref{t:main} for $G_0 = E_7(q)$. 

\subsection{$G_0 = E_6^{\e}(q)$, $q \geqs 3$}\label{s:e6}

Next we assume $G$ is an almost simple group with socle $G_0 = E_6^{\e}(q)$ and $H$ is a non-parabolic maximal subgroup. Throughout this section, we will assume that $q \geqs 3$; the case $q=2$ requires special attention and it will be handled separately in Section \ref{s:e62}. We will make extensive use of the information on conjugacy classes of semisimple and unipotent elements in \cite[Section 2]{FJ} and \cite[Table 22.2.3]{LSbook}.
In addition, we will also refer to \cite[Table 2]{CW}, which gives information on the centralizers of semisimple elements in $\bar{G}$ of order at most $7$. Throughout this section, we set $H_0 = H \cap G_0$.

\begin{lem}\label{l:e6_bd}
If $G_0 = E_6^{\e}(q)$, $q \geqs 3$ and $|H| \leqs q^{28}$, then $b(G) \leqs 5$.
\end{lem}

\begin{proof}
Define $\a,\b$ and $\gamma$ as in \eqref{e:alpha}. We begin by estimating $\a$, so let $x \in G$ be a semisimple element of prime order. As noted in the proof of \cite[Lemma 4.15]{BLS}, there are fewer than $2q^{32}$ such elements in $G$ with  
$\dim x^{\bar{G}} < 40$, and \cite[Theorem 2]{LLS} gives ${\rm fpr}(x) \leqs 2q^{-12}$. If $\dim x^{\bar{G}} \geqs 40$ then $|x^G|>\frac{1}{2}q^{40}$ and by applying Proposition \ref{p:bdd} and the given bound $|H| \leqs q^{28}$, we deduce that 
$$\a<2q^{32}(2q^{-12})^5+\frac{1}{2}q^{40}(2q^{-12})^5<q^{-17}.$$

Next let us turn to $\b$. As above, the contribution from unipotent elements $x \in G$ with $\dim x^{\bar{G}} \geqs 40$ is less than $\frac{1}{2}q^{40}(2q^{-12})^5$, so we may assume that $x$ belongs to one of the classes labelled $A_1$ or $A_1^2$. If $x$ is in $A_1$ then ${\rm fpr}(x) \leqs 2q^{-6}$ by \cite[Theorem 2]{LLS} and there are fewer than $2q^{22}$ such elements in $G$, so the contribution from long root elements is less than $2q^{22}(2q^{-6})^5 = 2^6q^{-8}$. Now assume $x$ is in the $A_1^2$ class, so $|x^G|>\frac{1}{2}q^{32}$. We claim that $H_0$ contains at most $q^{25}$ such elements, so their contribution is at most $\frac{1}{2}q^{32}(2q^{-7})^5 = 2^4q^{-3}$ by Proposition \ref{p:bdd}. This is obvious if $|H_0| \leqs q^{25}$, so let us assume $q^{25}<|H_0| \leqs q^{28}$. By arguing as in the proof of \cite[Lemma 4.13]{BLS}, which proceeds by inspecting the various possibilities for $H$ arising in \cite[Theorem 2]{LS90}, it is not difficult to rule out the existence of a maximal subgroup $H$ with $|H_0| \leqs q^{28}$ and $i_p(H_0)>q^{25}$. This establishes the claim and we conclude that
$$\b<\frac{1}{2}q^{40}(2q^{-12})^5 + 2^6q^{-8} + 2^4q^{-3} < 17q^{-3}.$$

To complete the proof, it remains to estimate the contribution $\gamma$ from field, graph and graph-field automorphisms. By Proposition \ref{p:bdd}, the contribution from elements $x \in G$ with $|x^G|>\frac{1}{6}q^{39}$ is less than $\frac{1}{6}q^{39}(6q^{-11})^5 = 6^4q^{-16}$, so we may assume $|x^G| \leqs \frac{1}{6}q^{39}$, which implies that $x$ is an involutory graph automorphism with $C_{\bar{G}}(x) = F_4$. We claim that 
$${\rm fpr}(x) \leqs \frac{1}{q^6-q^3+1},$$ 
so the contribution from these elements is at most $2q^{26}(q^6-q^3+1)^{-5} <q^{-3}$. For $\e=-$, the claim follows from \cite[Theorem 2]{LLS}, so we may assume $\e=+$ and it suffices to show that $i_2(H \setminus H_0) \leqs q^{19}$. Following the proof of \cite[Lemma 4.13]{BLS}, it is straightforward to reduce to the case where $H$ is an almost simple group of Lie type in characteristic $p$, and we can easily rule out the existence of such a subgroup with $|H| \leqs q^{28}$ and $i_2(H \setminus H_0) > q^{19}$. Therefore,
$\gamma < 6^4q^{-16} + q^{-3} < 2q^{-3}$ and we conclude that
$$\what{Q}(G,5)< 19q^{-3} + q^{-17}.$$
This completes the proof of the lemma.
\end{proof}

Next we determine the cases that arise with $|H|>q^{28}$ (recall that we are assuming $H$ is non-parabolic throughout this section).

\begin{lem}\label{l:e6q28}
If $G_0 = E_6^{\e}(q)$, $q \geqs 3$ and $|H|>q^{28}$, then one of the following holds:
\begin{itemize}\addtolength{\itemsep}{0.2\baselineskip}
\item[{\rm (i)}] $H = N_G(\bar{H}_{\s})$ with $\bar{H}^0 = F_4$, $T_1D_5$, $A_1A_5$, $T_2D_4$ or $C_4$ ($p \ne 2$);
\item[{\rm (ii)}] $\e=+$, $q=q_0^2$ and $H_0 = E_6^{\pm}(q_0)$.
\end{itemize}
\end{lem}

\begin{proof}
If $|H|>q^{32}$ then the result follows from \cite[Lemma 4.14]{BLS}, so we may assume $q^{28} < |H| \leqs q^{32}$. Since $|H|^3>|G|$, we can proceed as in the proof of \cite[Lemma 5.3]{AB} to rule out any additional subgroups. In particular, note that if $q=q_0^3$ and $H$ is a subfield subgroup of type $E_6^{\e}(q_0)$ then 
$$|H| \leqs 2(3,q_0-\e)\log_2q.|E_6^{\e}(q_0)|< q^{28}.$$
The result follows.
\end{proof}

Some of the cases arising in Lemma \ref{l:e6q28} have already been handled in \cite{BLS}. Indeed, the proof of \cite[Proposition 4.21]{BLS} shows that $b(G) \leqs 5$ in (ii), and the same bound holds in (i) when $\bar{H}^0 = A_1A_5$ or $C_4$ (see \cite[Lemmas 4.16 and 4.20]{BLS}). Therefore, we may assume $H = N_G(\bar{H}_{\s})$ with $\bar{H} = F_4$, $T_1D_5$ or $T_2D_4.S_3$. We start with the case $\bar{H} = F_4$, which requires the most work.

\begin{lem}\label{l:e6q3f4}
If $G_0 = E_6^{\e}(q)$, $q \geqs 3$ and $H = N_G(\bar{H}_{\s})$ with $\bar{H}=F_4$, then $b(G) \leqs 5$.
\end{lem}

\begin{proof}
First observe that $\bar{H} = C_{\bar{G}}(\tau)$, where $\tau$ is an involutory graph automorphism of $\bar{G}$. Write $\what{Q}(G,5) = \a+\b+\gamma$ as before (see \eqref{e:alpha}) and note that 
$$|\bar{G}_{\s}:\bar{H}_{\s}| = q^{12}(q^5-\e)(q^9-\e).$$

\vs

\noindent \emph{Case 1. Semisimple elements.} We start by estimating $\a$. Let $x \in G$ be a semisimple element of prime order $r$ and note that $\tau$ induces an automorphism on $\bar{D} = C_{\bar{G}}(x)$. First assume $\dim x^{\bar{G}} \leqs 52$. Here \cite[Theorem 2]{LLS2} gives ${\rm fpr}(x) \leqs 1/q^6(q^6-q^3+1)$ and by inspecting  \cite{FJ} we calculate that $G$ contains fewer than $q^{56}$ such elements, hence their combined contribution to $\a$ is less than $q^{-3}$. To obtain a sufficient bound on the remaining contribution to $\a$, we need to improve the fixed point ratio estimates in \cite{LLS2}. To do this, we will work closely with the proof of \cite[Lemma 5.4]{LLS2}.

Suppose $\bar{D}^0$ has an $A_3$ factor. By the proof of  \cite[Lemma 5.4]{LLS2} we have  
$\bar{D}^0 = A_3T_3$ (otherwise ${\rm fpr}(x)=0$) and $C_{\bar{D}}(\tau) = C_2T_2$, so  
$${\rm fpr}(x) \leqs \frac{|C_{\bar{G}_{\s}}(x): C_{\bar{H}_{\s}}(x)|}{|\bar{G}_{\s}:\bar{H}_{\s}|}$$
with $|C_{\bar{H}_{\s}}(x)| \geqs |{\rm Sp}_{4}(q)|(q-1)^2$ and $|C_{\bar{G}_{\s}}(x)| \leqs |{\rm SU}_{4}(q)|(q+1)^3$, hence ${\rm fpr}(x) < q^{-18}$. By inspecting \cite{FJ}, one checks that 
$G$ contains fewer than $|E_6^{\e}(q): {\rm SL}_{4}(q)| < 2q^{63}$ 
such elements, so their contribution to $\a$ is less than $2q^{63}(q^{-14})^5 = 2q^{-7}$. Similarly, if $\bar{D}^0 = A_2^3$ then $r=3$ (so $q \geqs 4$), $|x^G|>\frac{1}{6}q^{54}$ and the proof of \cite[Lemma 5.4]{LLS2} gives $C_{\bar{H}}(x) = A_2\tilde{A}_2$. Since $H$ contains at most $4q^{36}$ such elements, Proposition \ref{p:bdd} implies that their contribution is less than $\frac{1}{6}q^{54}(24q^{-18})^5<q^{-25}$. 

By inspecting \cite{FJ}, the remaining possibilities for $\bar{D}^0$ are as follows:
$$\begin{array}{l|ccccccccc}
\bar{D}^0 & A_2^2A_1T_1 & A_2^2T_2 & A_2A_1^2T_2 & A_2A_1T_3 & A_1^3T_3 & A_2T_4 & A_1^2T_4 & A_1T_5 & T_6 \\ \hline
\dim x^{\bar{G}} & 58 & 60 & 62 & 64 & 66 & 66 & 68 & 70 & 72 \\
\end{array}$$
(Notice that $\bar{D}^0 \ne A_1^4T_2$ since $x$ has prime order; see \cite[14.1]{GL}.) As noted in the proof of \cite[Lemma 5.4]{LLS2}, there are at most 
$|W(E_6):W(F_4)|=45$ distinct $\bar{H}$-classes in $x^{\bar{G}}\cap \bar{H}$. Therefore, if we choose $x$ with $|C_{\bar{H}_{\s}}(x)|$ minimal, then
\begin{equation}\label{e:45}
{\rm fpr}(x) \leqs 45 \cdot \frac{|C_{\bar{G}_{\s}}(x): C_{\bar{H}_{\s}}(x)|}{|\bar{G}_{\s}:\bar{H}_{\s}|}  = 45 \cdot \frac{|\bar{D}_{\s} : C_{\bar{D}_{\s}}(\tau)|}{|\bar{G}_{\s}:\bar{H}_{\s}|}.
\end{equation}
We will use this upper bound to show that
$${\rm fpr}(x) < \left\{\begin{array}{ll}
q^{-16} & \mbox{if $\bar{D}^0 = A_1^2T_4$, $A_1T_5$ or $T_6$} \\
q^{-15} & \mbox{if $\bar{D}^0=A_1^3T_3$ or $A_2T_4$} \\
q^{-14} & \mbox{if $\bar{D}^0= A_2A_1^2T_2$ or $A_2A_1T_3$} \\
q^{-13} & \mbox{if $\bar{D}^0= A_2^2A_1T_1$ or $A_2^2T_2$} 
\end{array}\right.$$
for the remaining semisimple elements $x \in G$ that we are interested in. 

First we show that ${\rm fpr}(x) < q^{-16}$ if $\bar{D}^0 = A_1^2T_4$, $A_1T_5$ or $T_6$, so the contribution from these elements is less than $q^{78}(q^{-16})^5 = q^{-2}$. To justify the claim, first note that $r \geqs 11$ (see \cite[Table 2]{CW}, for example). Since $|C_{\bar{H}_{\s}}(x)| \geqs (q-1)^4$ and $|C_{\bar{G}_{\s}}(x)| \leqs |{\rm SL}_{2}(q^2)|(q+1)^4$, the bound in \eqref{e:45} is sufficient if $q \geqs 5$. Suppose $q=4$. Since $r$ has to divide the order of the centre of $C_{\bar{H}_{\s}}(x)$, we deduce that $|C_{\bar{H}_{\s}}(x)| \geqs (q^2+1)(q-1)^2$ (see \cite[Table III]{Shinoda}) and the result follows as before. Similarly, if $q=3$ then the information on semisimple centralizers in \cite[Table 9]{Shoji} and \cite{FJ} implies that $|C_{\bar{H}_{\s}}(x)| \geqs (q^3-1)(q-1)$, $|C_{\bar{G}_{\s}}(x)| \leqs |{\rm SL}_{2}(q^2)|(q^2-1)(q^2+q+1)$ and once again the desired bound follows via \eqref{e:45}.

Similar reasoning shows that ${\rm fpr}(x) < q^{-13}$ if $\bar{D}^0 = A_2^2A_1T_1$ or $A_2^2T_2$. For instance, if $\bar{D}^0 = A_2^2A_1T_1$ then by arguing as in the proof of \cite[Lemma 5.4]{LLS2} we deduce that 
$$|C_{\bar{H}_{\s}}(x)| \geqs |{\rm SL}_{3}(q)||{\rm SL}_{2}(q)|(q-1), \;\; |C_{\bar{G}_{\s}}(x)| \leqs |{\rm SU}_{3}(q^2)||{\rm SL}_{2}(q)|(q+1)$$ 
and the result follows from \eqref{e:45}. By carefully inspecting \cite{FJ}, one checks that there are at most $q^{63}$ such elements in $G$, so the contribution to $\a$ from the elements with $\bar{D}^0 = A_2^2A_1T_1$ or $A_2^2T_2$ is less than $q^{63}(q^{-13})^5 = q^{-2}$. 

Next we show that ${\rm fpr}(x)<q^{-14}$ if $\bar{D}^0 = A_2A_1^2T_2$ or $A_2A_1T_3$. Suppose $\bar{D}^0 = A_2A_1^2T_2$. By inspecting \cite{FJ}, we deduce that $q \geqs 4$ (that is, if $q=3$ then there are no prime order elements $x \in G$ with $\bar{D}^0 = A_2A_1^2T_2$). Now $C_{\bar{H}}(x) = C_{\bar{D}}(\tau) = A_2A_1T_1$ or $A_1^2T_2$, so
$$|C_{\bar{H}_{\s}}(x)| \geqs |{\rm SL}_{2}(q)|^2(q-1)^2, \;\; |C_{\bar{G}_{\s}}(x)| \leqs |{\rm SU}_{3}(q)||{\rm SL}_{2}(q^2)|(q+1)^2$$ 
and the desired bound holds. Similarly, if $\bar{D}^0 = A_2A_1T_3$ then $r \geqs 7$ and the bounds 
$$|C_{\bar{H}_{\s}}(x)| \geqs |{\rm SL}_{2}(q)|(q-1)^3, \;\; |C_{\bar{G}_{\s}}(x)| \leqs |{\rm SU}_{3}(q)||{\rm SL}_{2}(q)|(q+1)^3$$ 
are sufficient (via \eqref{e:45}) if $q \geqs 4$. Finally, if $q=3$ then we can replace the terms $(q-1)^3$ and $(q+1)^3$ in the above bounds by $(q^3-1)$ and $(q^3+1)$, respectively (since $r \geqs 7$), and the result follows as before. One checks that there are at most $q^{68}$ elements in $G$ with 
$\bar{D}^0 = A_2A_1^2T_2$ or $A_2A_1T_3$, so the contribution to $\a$ is less than $q^{68}(q^{-14})^5 = q^{-2}$.

Finally, let us assume $\bar{D}^0 = A_1^3T_3$ or $A_2T_4$. We claim that ${\rm fpr}(x) < q^{-15}$. In both cases, $r \geqs 7$ and
$$|C_{\bar{H}_{\s}}(x)| \geqs |{\rm SL}_{2}(q)|(q-1)^3, \;\; |C_{\bar{G}_{\s}}(x)| \leqs |{\rm SU}_{3}(q)|(q+1)^4,$$ 
so the claim follows from \eqref{e:45} if $q \geqs 4$. For $q=3$, we can replace the $(q-1)^3$ term in the lower bound for $|C_{\bar{H}_{\s}}(x)|$ by $(q^3-1)$ and the result follows. By inspecting \cite{FJ}, one checks that $G$ contains fewer then $q^{71}$ such elements, so their contribution is less than $q^{71}(q^{-15})^5 = q^{-4}$.

Bringing together all of the above estimates, we conclude that
\begin{equation}\label{e:al}
\a< 3q^{-2} +q^{-3} + q^{-4} + 2q^{-7}+q^{-25} < 4q^{-2}.
\end{equation}

\vs

\noindent \emph{Case 2. Unipotent elements.} Now let us turn to the contribution from unipotent elements. Note that the unipotent classes in $\bar{H}$, together with the corresponding classes in $\bar{G}$, are listed in \cite[Table A]{Law}. We refer the reader to \cite[Tables 22.2.3 and 22.2.4]{LSbook} for information on the conjugacy classes and centralizers of unipotent elements in $\bar{G}_{\s}$ and $\bar{H}_{\s}$.

First assume $p=2$, so $q \geqs 4$. If $x \in G$ is a long root element, then \cite[Theorem 2]{LLS2} gives ${\rm fpr}(x) \leqs (q^6-q^3+1)^{-1}$ and it is easy to check that there are fewer than $2q^{22}$ such elements in $G$. Similarly, if $x$ is in the $A_1^2$ class then $|x^G|>\frac{1}{2}q^{32}$ and there are at most $2q^{22}$ such elements in $H$, and for $x$ in the $A_1^3$ class we have $|x^G|>\frac{1}{2}q^{40}$ and $H$ contains fewer than $2q^{28}$ such elements. By applying Proposition \ref{p:bdd}, we conclude that 
$$\b<2q^{22}(q^6-q^3+1)^{-5} + \frac{1}{2}q^{32}(4q^{-10})^5 + \frac{1}{2}q^{40}(4q^{-12})^5<q^{-7}$$
when $p=2$.

Now assume $p$ is odd. First we claim that the contribution from unipotent elements with $\dim x^{\bar{G}} \geqs 54$ is at most $q^{-3}$. By inspecting \cite[Table A]{Law}, we see that $\delta(x) \geqs 18$. For $q \geqs 5$, the bound in \cite[(4.16)]{BLS} implies that ${\rm fpr}(x) < q^{-15}$ and so the claim follows from the fact that $G$ contains fewer than $q^{72}$ such elements. Now assume $q=3$, so $x$ has order $3$ and the bound on $\dim x^{\bar{G}}$ implies that $x$ is in the class $A_2^2A_1$ (for example, this follows by considering the Jordan form of $x$ on the Lie algebra of $\bar{G}$, as given in \cite[Table 6]{Law}). Here $\delta(x) = 18$ and $C_{\bar{G}}(x)$ is connected, so ${\rm fpr}(x) < q^{-13}$ by \cite[(4.16)]{BLS} and we note that there are at most $2q^{54}$ of these elements in $G$, so the contribution is less than $2q^{54}(q^{-13})^5$ and the claim follows. 

As noted in the proof of \cite[Lemma 4.19]{BLS}, the contribution from unipotent elements  with $48 \leqs \dim x^{\bar{G}} < 54$ is less than $2q^{52}(q^{-14})^5 = 2q^{-18}$. Finally, let us assume $\dim x^{\bar{G}}<48$ and ${\rm fpr}(x)>0$, so $x$ is in one of the $\bar{G}$-classes labelled $A_1$, $A_1^2$, $A_1^3$ or $A_2$ (see \cite[Table A]{Law}). As above, the contribution from long root elements is less than $2q^{22}(q^6-q^3+1)^{-5}$, with the remaining unipotent elements contributing at most $\sum_{i=1}^{3}a_i(b_i/a_i)^5$, where
$$a_1 = \frac{1}{2}q^{32}, \; a_2 = \frac{1}{2}q^{40}, \; a_3 = \frac{1}{4}q^{42},\; b_1 = 2q^{22}, \; b_2 = 2q^{28}, \; b_3 = 2q^{30}.$$
Therefore, 
\begin{equation}\label{e:be}
\b < q^{-3} + 2q^{-18} + 2q^{22}(q^6-q^3+1)^{-5} + \sum_{i=1}^{3}a_i(b_i/a_i)^5<q^{-2}
\end{equation}
if $p$ is odd. 

\vs

\noindent \emph{Case 3. Remaining elements.} It remains to estimate the contribution 
$\gamma$ from field, graph and graph-field automorphisms. As in the proof of \cite[Lemma 4.19]{BLS}, the contribution from involutory graph and graph-field automorphisms is less than $4q^{39}(q^{-12})^5 = 4q^{-21}$. Similarly, the contribution from odd order field automorphisms is at most 
$$\log_2q.q^{78}(12q^{-52/3})^5 < q^{-2}$$
(the inequality holds since we may assume $q \geqs 8$) and that from involutory graph automorphisms is less than 
$$2q^{42}(q^{-11})^5 + 2q^{26}(12q^{-10})^5.$$
Therefore
\begin{equation}\label{e:ga}
\gamma< 4q^{-21} + q^{-2} + 2q^{42}(q^{-11})^5 + 2q^{26}(12q^{-10})^5 < q^{-1}
\end{equation}
and by bringing together the bounds in \eqref{e:al}, \eqref{e:be} and \eqref{e:ga}, we conclude that
$$\what{Q}(G,5) < q^{-1} + 5q^{-2}.$$
The result follows.
\end{proof}

\begin{lem}\label{l:e6q3t1d5}
If $G_0 = E_6^{\e}(q)$, $q \geqs 3$ and $H = N_G(\bar{H}_{\s})$ with $\bar{H}=T_1D_5$, then $b(G) \leqs 5$.
\end{lem}

\begin{proof}
We proceed as in the proof of the previous lemma, estimating $\a$, $\b$ and $\gamma$ in turn. The argument closely follows the proof of \cite[Lemma 4.17]{BLS}. 

Let $x \in G$ be a semisimple element of prime order $r$ and set $\bar{D} = C_{\bar{G}}(x)$, $z = \dim Z(\bar{D}^0)$ and $\delta(x) = \dim x^{\bar{G}} - \dim (x^{\bar{G}} \cap \bar{H})$. As noted in the proof of the previous lemma, there are fewer than $q^{56}$ semisimple elements with $\dim x^{\bar{G}} \leqs 52$, so \cite[Theorem 2]{LLS2} implies that their contribution is less than $q^{56}(2q^{-12})^5 = 2^5q^{-4}$. 

Next we claim that ${\rm fpr}(x) < q^{-16}$ if $\dim x^{\bar{G}} \geqs 60$. For elements with $z \leqs 4$, this follows immediately from the upper bound in \cite[(4.11)]{BLS}, since $\delta(x) \geqs 24$ (the latter bound is established in the proof of \cite[Lemma 4.17]{BLS}). On the other hand, if $z>4$ then $\bar{D}^0 = A_1T_5$ or $T_6$, so 
$$|x^{G}|>\frac{1}{2}\left(\frac{q}{q+1}\right)^5q^{70}$$
and the desired bound holds since $|H \cap \bar{G}_{\s}| < q^{47}$. Therefore, the contribution from these elements is less than $q^{78}(q^{-16})^5 = q^{-2}$. 

To complete the analysis of semisimple elements, we may assume that $54 \leqs \dim x^{\bar{G}} \leqs 58$, in which case $\bar{D}^0 = A_2^3$, $A_3A_1T_2$ or $A_2^2A_1T_1$ (note that $\bar{D}^0 \ne A_3A_1^2T_1$ by \cite[14.1]{GL}). First observe that $q \geqs 4$. This is clear if $\bar{D}^0 = A_2^3$ since $r=3$; in the other two cases, $r \geqs 5$ must divide the order of the centre of $C_{\bar{G}_{\s}}(x)$ and by inspecting \cite{FJ} we deduce that $q \geqs 4$. In addition, $z \leqs 2$ and $\delta(x) \geqs 20$ by \cite[Theorem 2]{LLS}, so \cite[(4.11)]{BLS} implies that ${\rm fpr}(x) < q^{-13}$. Therefore, the combined contribution from these elements is less than $q^{62}(q^{-13})^5 = q^{-3}$. We conclude that
$$\a< 2^5q^{-4} + q^{-2} + q^{-3}< 2q^{-2}.$$

Next let us consider the contribution from unipotent elements. First assume $p=2$. As explained in the proof of \cite[Lemma 4.17]{BLS}, we can determine the $\bar{G}$-class of each involution in $\bar{H}$ by considering the restriction $V_{27}{\downarrow}D_5$, where $V_{27}$ is one of the minimal modules for $E_6$. In terms of the standard Aschbacher-Seitz \cite{AS} notation for involutions in $D_5$, we find that $a_2 \in A_1$, $c_2,a_4 \in A_1^2$ and $c_4 \in A_1^3$. This implies that $\b < \sum_{i=1}^{3}u_i(v_i/u_i)^5$, where
\begin{equation}\label{e:ui}
u_1 = \frac{1}{2}q^{22}, \; u_2 = \frac{1}{2}q^{32}, \; u_3 = \frac{1}{2}q^{40}, \; v_1 = 2q^{14},\; v_2 = 2q^{16}(q^4+1),\; v_3 = 2q^{24}.
\end{equation}
In particular, $\b<q^{-13}$.

Now assume $p$ is odd. Once again, as noted in the proof of \cite[Lemma 4.17]{BLS}, we can determine the $\bar{G}$-class of each unipotent element in $\bar{H}$ and it is routine to check   that the upper bound on ${\rm fpr}(x)$ presented in \cite[(4.10)]{BLS} is satisfied. If $\dim x^{\bar{G}} \geqs 50$ then $\delta(x) \geqs 20$ and we deduce that ${\rm fpr}(x) < q^{-16}$ if $q \geqs 5$. On the other hand, if $q=3$ then $r=3$ and we see that $x \in A_2A_1^2$ is the only possibility, given the condition on $\dim x^{\bar{G}}$. Here \cite[(4.10)]{BLS} gives ${\rm fpr}(x) < q^{-14}$ and there are fewer than $2q^{50}$ such elements in $G$. Therefore, for any $q$, it follows that the contribution to $\b$ from these unipotent elements is less than $q^{72}(q^{-16})^5 = q^{-8}$. Finally, suppose $\dim x^{\bar{G}} < 50$, so $x$ belongs to one of the classes labelled $A_1$, $A_1^2$, $A_2$, $A_1^3$ or $A_2A_1$ (no unipotent elements in $\bar{H}$ are in the $\bar{G}$-class $A_2^2$). By applying the bounds in the proof of \cite[Lemma 4.17]{BLS}, we conclude that
$$\b < q^{-8} + 2q^{22}(2q^{-6})^5 + 3q^{32}(q^{-9})^5 + 3q^{46}(q^{-10})^5 < q^{-2}.$$

Finally, let us consider $\gamma$. By directly applying the bounds presented in the proof of \cite[Lemma 4.17]{BLS}, we deduce that
$$\gamma < 4q^{39}(q^{-12})^5 + \log_2q.q^{78}(q^{-19})^5 + 2q^{42}(q^{-11})^5 + 2q^{26}(4(q+1)q^{-17})^5 < q^{-12}.$$

Bringing together the above bounds, we conclude that
$\what{Q}(G,5) < 3q^{-2} + q^{-12}$
and the result follows.
\end{proof}

\begin{lem}\label{l:e6q3t2d4}
If $G_0 = E_6^{\e}(q)$, $q \geqs 3$ and $H = N_G(\bar{H}_{\s})$ with $\bar{H}=T_2D_4.S_3$, then $b(G) \leqs 5$.
\end{lem}

\begin{proof}
This is similar to the proof of the previous lemma, working with the proof of \cite[Lemma 4.18]{BLS} in place of \cite[Lemma 4.17]{BLS}. Note that one of ${\rm P\O}_{8}^{+}(q)$ or ${}^3D_4(q)$ is a composition factor of $\bar{H}_{\s}$ (see \cite[Table 5.1]{LSS}).

First let $x \in G$ be a semisimple element of prime order $r$ and define $\bar{D}$, $z$ and $\delta(x)$ as before. We claim that ${\rm fpr}(x)<q^{-16}$ if $\dim x^{\bar{G}} \geqs 48$. For $q \geqs 4$, this follows immediately from the upper bound in \cite[(4.13)]{BLS} since $z \leqs 6$ and $\delta(x) \geqs 26$ by \cite[Theorem 2]{LLS2}. The same bound is also sufficient if $q=3$ and $z \leqs 3$. Finally, if $q=3$ and $z \geqs 4$ then 
$$|x^G|>\frac{1}{2}\left(\frac{q}{q+1}\right)^4q^{66}$$
and the claim follows since $|H \cap \bar{G}_{\s}| < q^{32}$. This justifies the claim and it follows that the contribution to $\a$ from these elements is less than $q^{78}(q^{-16})^5 = q^{-2}$. By inspecting \cite{FJ}, one checks that there are fewer than $q^{44}$ semisimple elements in $G$ with $\dim x^{\bar{G}}<48$ and \cite[Theorem 2]{LLS2} gives ${\rm fpr}(x) \leqs 2q^{-12}$. Therefore,
$$\a < q^{-2} + q^{44}(2q^{-12})^5 < 2q^{-2}.$$

Next let us consider unipotent elements. First assume $p=2$. If ${}^3D_4(q)$ is a composition factor of $\bar{H}_{\s}$, then the proof of \cite[Lemma 4.18]{BLS} implies that 
$$\b <  2q^{22}(2q^{-12})^5 + 2q^{40}(4q^{-24})^5.$$
Similarly, if $\O_{8}^{+}(q)$ is a composition factor, then $\b < \sum_{i=1}^{3}u_i(v_i/u_i)^5$, where
$$u_1 = \frac{1}{2}q^{22}, \; u_2 = \frac{1}{2}q^{32}, \; u_3 = \frac{1}{2}q^{40},\; v_1 = 2q^{10}, \; v_2 = 2q^7(2q^5+1),\; v_3 = 2q^{15}(q+1).$$
In both cases, these estimates imply that $\b<q^{-33}$ when $p=2$.

Now assume $p$ is odd. As explained in the proof of \cite[Lemma 4.18]{BLS}, if $\dim x^{\G} \geqs 40$ then ${\rm fpr}(x) < q^{-16}$ and so the contribution from these elements is less than $q^{72}(q^{-16})^5 = q^{-8}$ and it remains to consider the elements in the classes $A_1$ and $A_1^2$. If $x$ is in $A_1^2$ then the proof of \cite[Lemma 4.18]{BLS} gives ${\rm fpr}(x) < q^{-12}$, whereas \cite[Theorem 2]{LLS2} gives ${\rm fpr}(x) \leqs 2q^{-6}$ for a long root element $x$. Therefore, 
$$\b < q^{-8} + 2q^{32}(q^{-12})^5 + 2q^{22}(2q^{-6})^5 < q^{-4}$$
when $p$ is odd.

Finally, let us consider $\gamma$. From the bounds given in the proof of \cite[Lemma 4.18]{BLS}, we deduce that the contribution from field and graph-field automorphisms is less than 
$$4q^{39}(q^{-12})^5 + \log_2q.q^{78}(q^{-20})^5.$$
Now assume $x$ is an involutory graph automorphism. If $C_{\bar{G}}(x) \ne F_4$ then the proof of \cite[Lemma 4.18]{BLS} gives ${\rm fpr}(x) < q^{-17}$, so the contribution to $\gamma$ is less than $2q^{42}(q^{-17})^5$. Finally, suppose $C_{\bar{G}}(x) = F_4$, so
$$|x^G| = \frac{q^{12}(q^5-\e)(q^9-\e)}{(3,q-\e)}.$$
Here the proof of \cite[Lemma 4.18]{BLS} gives ${\rm fpr}(x) < q^{-5}$, which is not sufficient. We claim that ${\rm fpr}(x) < q^{-6}$, which implies that the contribution from these elements is less than $2q^{26}(q^{-6})^5 = 2q^{-4}$. To justify the claim, first observe that 
\begin{equation}\label{e:xgh}
|x^G \cap H| \leqs i_2(\bar{H}_{\s}) \leqs 2(q+1)^2\cdot i_2({\rm Aut}(S)),
\end{equation}
where $S = {\rm P\O}_{8}^{+}(q)$ or ${}^3D_4(q)$ according to the structure of $\bar{H}_{\s}$. As noted in \cite[(4.8)]{BLS}, this gives $|x^G \cap H| < 4(q+1)^3q^{15}$ and one checks that this bound is sufficient if $q \geqs 5$. For $q \in \{3,4\}$, we compute $i_2({\rm Aut}(S))$ precisely and then use the upper bound in \eqref{e:xgh}. For example, if $S = {}^3D_4(4)$ then 
$$i_2({\rm Aut}(S)) = i_2({}^3D_4(4)) + |{}^3D_4(4) : {}^3D_4(2)| = 4632563455$$
and the claim follows. We conclude that
$$\gamma< 4q^{39}(q^{-12})^5 + \log_2q.q^{78}(q^{-20})^5 + 2q^{42}(q^{-17})^5 + 2q^{-4} < 3q^{-4}.$$

Therefore, $\what{Q}(G,5) < 2q^{-2} + 4q^{-4}$ and the proof of the lemma is complete.
\end{proof}

This completes the proof of Theorem \ref{t:main} for $G_0 = E_6^{\e}(q)$ with $q \geqs 3$. 

\subsection{$G_0 = E_6^{\e}(2)$}\label{s:e62}

Here we handle the special cases with $G_0 = E_6(2)$ or ${}^2E_6(2)$. We refer the reader to \cite[Table 9]{BLS} for useful information on the classes of elements in $\bar{G}_{\s}$ of odd prime order. The maximal subgroups of $G$ are determined in \cite{KW} and \cite{W18} for $G_0 = E_6(2)$ and ${}^2E_6(2)$, respectively (in the latter case, it turns out that the list of maximal subgroups presented in the Atlas \cite{Atlas} is complete).

We begin by extending Lemma \ref{l:e6_bd} to the case $q=2$.

\begin{lem}\label{l:e6_bd2}
If $G_0 = E_6^{\e}(2)$ and $|H| \leqs 2^{28}$, then $b(G) \leqs 5$.
\end{lem}

\begin{proof}
From the list of maximal subgroups in \cite{Atlas, KW}, we first observe that the condition $|H| \leqs 2^{28}$ implies that $|H|< 2^{27}$. Write $\what{Q}(G,5) = \a+\b+\gamma$ as before.

We start by estimating $\a$. Let $x \in G$ be an element of odd prime order $r$. Clearly, the contribution from the elements with $|x^G|>2^{45}$ is less than $2^{45}(2^{27}/2^{45})^5 = 2^{-45}$. Now assume $|x^G| \leqs 2^{45}$, so $r=3$. If $\e=+$, or if $\e=-$ and $G = G_0$ or $G_0.2$, then $|x^G|>2^{41}$ and thus the contribution to $\a$ is less than $2^{-29}$. If $\e=-$ and $G = G_0.3$ or $G_0.S_3$ then $|x^G|>2^{31}$ and by inspecting the possibilities for $H$ in \cite{Atlas}, one checks that $i_3(H) < 2^{24}$. For example, if $H_0 = {\rm U}_{3}(8).3$ then 
$$i_3(H) \leqs i_3({\rm Aut}({\rm U}_{3}(8)) \times 3) = 1266434 < 2^{21}.$$
This means that the contribution from these elements is less than $2^{31}(2^{24}/2^{31})^5 = 2^{-4}$ and we conclude that $\a< 2^{-4}+2^{-45}$.

Next let us turn to the contribution from unipotent involutions. There are three $\bar{G}$-classes to consider, labelled $A_1$, $A_1^2$ and $A_1^3$. If $x \in A_1$ then ${\rm fpr}(x) \leqs 2^{-5}$ by \cite[Theorem 2]{LLS} and there are fewer than $2^{23}$ such elements in $G$, so their contribution to $\b$ is less than $2^{23}(2^{-5})^5 = 2^{-2}$. Now assume $x$ is in $A_1^2$ or $A_1^3$. If $\e=+$ then $|x^G|>2^{33}$ and $|H_0|<2^{26}$, so these elements contribute less than $2^{33}(2^{26}/2^{33})^5 = 2^{-2}$. 
Similarly, if $\e=-$ then $|x^G|>2^{31}$, $|H_0|<2^{24}$ and the contribution is at most $2^{-4}$. We conclude that $\b<2^{-1}$.

Finally, let us assume $x \in G$ is an involutory graph automorphism. If $C_{\bar{G}}(x) \ne F_4$ then $|x^G|>2^{40}$ and so the contribution to $\gamma$ from these elements is less than $2^{40}(2^{27}/2^{40})^5 = 2^{-25}$. On the other hand, if $C_{\bar{G}}(x) = F_4$ then $|x^G|<2^{27}$ and by combining \cite[Theorem 2]{LLS2} and \cite[Lemma 4.12]{BLS} we see that 
${\rm fpr}(x)\leqs 1/57$, so the contribution here is at most $2^{27}(1/57)^5<2^{-2}$. This gives $\gamma<2^{-2}+2^{-45}$ and we conclude that 
$$\what{Q}(G,5)< 2^{-1} + 2^{-2}+2^{-4}+2^{-44}<1$$
as required.
\end{proof}

\begin{lem}\label{l:e6_bd20}
If $G_0 = E_6(2)$ and $|H| > 2^{28}$, then $b(G) \leqs 5$.
\end{lem}

\begin{proof}
By inspecting \cite[Table 1]{KW}, we deduce that $H_0 = F_4(2)$, $S_3 \times {\rm L}_{6}(2)$ or $(7 \times {}^3D_4(2)).3$. In fact, in view of \cite[Lemmas 4.16 and 4.19]{BLS}, we only need to handle the latter case. Here $H = (D_{14} \times {}^3D_4(2)).3$ if $G = G_0.2$, and it will be useful to observe that $|H_0|<2^{33}$ and $i_3(H) \leqs i_3(7.3 \times {}^3D_4(2).3) < 2^{26}$. Write $\what{Q}(G,5) = \a+\b+\gamma$ as usual.

First we estimate $\a$. The contribution from the elements with $|x^G|>2^{47}$ is less than $2^{47}(2^{33}/2^{47})^5 = 2^{-23}$. On the other hand, if $|x^G| \leqs 2^{47}$ then $x$ has order $3$ and $|x^G|>2^{41}$, so the contribution from these elements is less than $2^{41}(2^{26}/2^{41})^5 = 2^{-34}$. 

Finally, let us consider $\b+\gamma$, which is the contribution from involutions. The subgroup 
${}^3D_4(2)<H_0$ has two classes of involutions, labelled $A_1$ and $A_1^3$, and the labels stay the same when viewed as involutions in $G_0$ (see the proof of \cite[Lemma 4.18]{BLS}). If $x$ is in $A_1$ then $|x^G|>2^{22}$ and there are at most $2^{10}$ such elements in $H$, so the contribution is less than $2^{22}(2^{-12})^5 = 2^{-38}$. If $x \in G$ is any other involution, then $|x^G| \geqs 2^{12}(2^5-1)(2^9-1)=a$ and $i_2(H) = 556927=b$, so the contribution here is at most $a(b/a)^5<2^{-8}$. 

We conclude that $\what{Q}(G,5)<1$.
\end{proof}

\begin{lem}\label{l:e6_bd30}
If $G_0 = {}^2E_6(2)$ and $|H| > 2^{28}$, then $b(G) \leqs 5$.
\end{lem}

\begin{proof}
According to the Atlas \cite{Atlas}, $H_0$ is one of the following:
$${\rm (a)} \; F_4(2), \; {\rm (b)} \; {\rm Fi}_{22}, \; {\rm (c)} \; S_3 \times {\rm U}_{6}(2), \; {\rm (d)} \; \Omega_{7}(3),$$
$${\rm (e)} \; \O_{10}^{-}(2), \; {\rm (f)} \; (3 \times \O_{8}^{+}(2).3).2, \; {\rm (g)} \; {}^3D_4(2).3.$$
By inspecting \cite{BLS}, we immediately deduce that $b(G) \leqs 5$ in cases (a)--(d), so it remains to consider (e), (f) and (g). We define $\a$, $\b$ and $\gamma$ in the usual fashion.

\vs

\noindent \emph{Case 1. $H_0 = \O_{10}^{-}(2)$.} As explained in the proof of Lemma \ref{l:e6q3t1d5}, we have $\b<\sum_{i=1}^3u_i(v_i/u_i)^5$, where the $u_i$ and $v_i$ terms are defined in \eqref{e:ui}, whence $\b < 2^{-8}$.

Next consider $\a$ and let $x \in G$ be a semisimple element of prime order $r$. If $r>3$ then $|x^G|>2^{58}$ and so the contribution to $\a$ is less than $2^{58}(2^{45}/2^{58})^5 = 2^{-7}$ since $|H_0|<2^{45}$. Now assume $r=3$ and note that $i_3(H) \leqs i_3(\O_{10}^{-}(2) \times 3) < 2^{32}$. Therefore, the contribution from the elements with $|x^G|>2^{41}$ is less than $2^{41}(2^{32}/2^{41})^5 = 2^{-4}$. Now assume $|x^G| \leqs 2^{41}$, in which case $x$ is in one of the classes labelled $3D$ and $3E$ in \cite{Atlas} and there are fewer than $2^{33}$ such elements in $G$. Since $C_{\bar{G}}(x) = D_5T_1$, the proof of \cite[Lemma 4.7]{LLS2} gives ${\rm fpr}(x) \leqs 2^{-11}$ and we conclude that the contribution to $\a$ is less than $2^{33}(2^{-11})^5 = 2^{-22}$. Therefore, $\a< 2^{-4}+2^{-7}+2^{-22}$.

Finally, let $x \in G$ be an involutory graph automorphism. If $C_{\bar{G}}(x) \ne F_4$ then $|x^G|<2^{27}$ and \cite[Theorem 2]{LLS2} gives ${\rm fpr}(x) \leqs 1/57$, so the contribution from these elements is less than $2^{27}(1/57)^5<2^{-2}$. On the other hand, if $C_{\bar{G}}(x) \ne F_4$ then $|x^G|>2^{40}$ and $i_2(H) \leqs i_2(O_{10}^{-}(2)) <2^{26}$, which gives a contribution of at most $2^{40}(2^{26}/2^{40})^5=2^{-30}$.

We conclude that $\what{Q}(G,5) < 2^{-2}+2^{-4}+2^{-7}+2^{-8}+2^{-22}+2^{-30}<1$.

\vs

\noindent \emph{Case 2. $H_0 = (3 \times \O_{8}^{+}(2).3).2$ or ${}^3D_4(2).3$.} We handle both cases simultaneously. Suppose $x \in G$ has odd prime order $r$. If $r=3$ then $|x^G|>2^{31}=a$ and we calculate that 
$$i_3(H) \leqs i_3(3^2 \times \O_{8}^{+}(2).3) = 18060488 = b.$$
Therefore, the contribution to $\a$ from elements of order $3$ is less than $a(b/a)^5<2^{-3}$. On the other hand, if $r>3$ then $|x^G|>2^{58}$ and so the remaining contribution is less than $2^{58}(2^{32}/2^{58})^5 =  2^{-72}$ since $|H_0|<2^{32}$. It follows that $\a<2^{-3}+2^{-72}$.

Finally, suppose $x \in G$ is an involution. If $x$ is a long root element then $|x^G|>2^{21}$ and by applying \cite[Proposition 1.13(ii)]{LLS} we deduce that $|x^G \cap H|<2^{11}$, so the contribution from these elements is less than $2^{21}(2^{11}/2^{21})^5 = 2^{-29}$. As in the previous case, the contribution from graph automorphisms with $C_{\bar{G}}(x) = F_4$ is at most $2^{-2}$. For the remaining involutions we have $|x^G|>2^{31}$ and their contribution is less than $2^{31}(2^{20}/2^{31})^5 = 2^{-24}$ since $i_2(H) < 2^{20}$. Therefore, $\b+\gamma < 2^{-2}+2^{-24}+2^{-29}$ and the result follows.
\end{proof}

\subsection{$G_0 = F_4(q)$}\label{s:f4}

In this final section we complete the proof of Theorem \ref{t:main} by handling the groups with socle $G_0 = F_4(q)$. We adopt the same notation as before. In particular, we will write $\what{Q}(G,5) = \a+\b+\gamma$. We refer the reader to \cite{Shinoda,Shoji} for information on the conjugacy classes in $F_4(q)$ (also see \cite[Table 22.2.4]{LSbook} and \cite{Lu2} for useful data on unipotent and semisimple classes, respectively). The case $q=2$ requires special attention and it will be treated separately in Lemma \ref{l:f42}. 

\begin{lem}\label{l:f4_bd}
If $G_0 = F_4(q)$, $q \geqs 3$ and $|H_0| \leqs q^{17}$, then $b(G) \leqs 5$. 
\end{lem}

\begin{proof}
We proceed as in the proof of \cite[Lemma 4.22]{BLS}. Suppose $x \in G_0$ has prime order. If $\dim x^{\bar{G}} \geqs 28$ then $|x^G|>q^{28}$ and so the contribution from these elements is less than $q^{28}(q^{17}/q^{28})^5 = q^{-27}$. Now assume $\dim x^{\G}<28$. If $x$ is semisimple, then $C_{\bar{G}}(x) = B_4$, $p \ne 2$ and $x$ is an involution. There are fewer than $2q^{16}$ such elements in $G$ and \cite[Theorem 2]{LLS2} gives ${\rm fpr}(x) \leqs 2q^{-5}$, so the contribution to $\a$ is less than $2q^{16}(2q^{-5})^5 = 2^6q^{-9}$. 

Now assume $x$ is unipotent and $\dim x^{\G}<28$, so $x$ belongs to one of the $\bar{G}$-classes labelled $A_1$, $\tilde{A}_1$ or $(\tilde{A}_1)_2$ (the latter only if $p=2$). First assume $x$ is in $A_1$ (or $\tilde{A}_1$ if $p=2$). Here ${\rm fpr}(x) \leqs (q^4-q^2+1)^{-1}$ by \cite[Theorem 2]{LLS2} and there are fewer than $3q^{16}$ such elements in $G$, so the contribution to $\b$ is less than $3q^{16}(q^4-q^2+1)^{-5} < q^{-2}$. Now suppose $x$ is in one of the classes $\tilde{A}_1$ (with $p \ne 2$) or  $(\tilde{A}_1)_2$ ($p=2$). Then $|x^G| \geqs \frac{1}{2}q^3(q^3-1)(q^4+1)(q^{12}-1) = f(q)$ and thus the contribution to $\b$ is at most $f(q) \cdot (q^{17}/f(q))^5<2q^{-1}$. Therefore,
$$\a+\b < 2q^{-1} + q^{-2} + 2^6q^{-9} + q^{-27}.$$

Finally, let us assume $x \in G$ is a field and graph-field automorphism. If $x$ has order $2$ then ${\rm fpr}(x) \leqs q^{-6}$ by \cite[Theorem 2]{LLS2} and there are fewer than $2q^{26}$ such elements in $G$, so the contribution to $\gamma$ is at most $2q^{26}(q^{-6})^5 = 2q^{-4}$. Similarly, if $x$ is a field automorphism of odd prime order then $|x^G|>\frac{1}{2}q^{34}$ and we note that $|H| \leqs 2\log_2q.q^{17}$, so the contribution from these elements is less than
$$\frac{1}{2}q^{34}(4\log_2q.q^{-17})^5 < q^{-43}.$$
Therefore, 
$$\what{Q}(G,5)<2q^{-1} + q^{-2} + 2q^{-4}+ 2^6q^{-9} + q^{-27} +q^{-43}<1$$
and thus $b(G) \leqs 5$.
\end{proof}

\begin{lem}\label{l:f4q17}
If $G_0 = F_4(q)$, $q \geqs 3$ and $|H_0| > q^{17}$, then one of the following holds:
\begin{itemize}\addtolength{\itemsep}{0.2\baselineskip}
\item[{\rm (i)}] $H = N_G(\bar{H}_{\s})$, where $\bar{H}^0 = A_1C_3$, $B_4$, $C_4$ ($p=2$ only) or $D_4$;
\item[{\rm (ii)}] $H = N_G(\bar{H}_{\s})$, where $\bar{H}^0 = C_2^2$ and $p=2$;
\item[{\rm (iii)}] $H_0 = F_4(q_0)$, where $q=q_0^k$ and $k \in \{2,3\}$;
\item[{\rm (iv)}] $H_0 = {}^2F_4(q)$, where $q=2^{2a+1}$ and $a \geqs 1$;
\item[{\rm (v)}] $H = {}^3D_4(2).3$ and $q=3$.
\end{itemize}
\end{lem}

\begin{proof}
By \cite[Theorem 2]{LS90}, one of the following holds:
\begin{itemize}\addtolength{\itemsep}{0.2\baselineskip}
\item[{\rm (a)}] $H = N_G(\bar{H}_{\s})$ for some $\s$-stable closed subgroup $\bar{H}$ of positive dimension;
\item[{\rm (b)}] $H$ is an exotic local subgroup (as in \cite[Theorem 1]{CLSS});
\item[{\rm (c)}] $H_0 = {}^2F_4(q)$ or $F_4(q_0)$, where $q=q_0^k$ and $k$ is a prime;
\item[{\rm (d)}] $H$ is almost simple and not of type (a) or (c).
\end{itemize}
Note that the subgroups with $|H|>q^{22}$ are determined in \cite[Lemma 4.23]{BLS}, which gives the examples in parts (i), (iii) (with $k=2$) and (iv). 

The relevant maximal rank subgroups of the form $H = N_G(\bar{H}_{\s})$ are determined in \cite{LSS} and it is easy to check that the only additional examples with $|H_0| > q^{17}$ are those described in part (ii) of the lemma (here $p=2$ and $G$ must contain a graph-field automorphism, as a condition for maximality). If $H =N_G(\bar{H}_{\s})$ is not of maximal rank, then 
$$|H_0| \leqs |{\rm SL}_{2}(q)||G_2(q)| < q^{17},$$
so none of these subgroups arise. The main theorem of \cite{CLSS} implies that there are no exotic local subgroups with $|H_0|>q^{17}$ and it is easy to see that $|F_4(q_0)|>q^{17}$ if and only if $q=q_0^k$ with $k =2$ or $3$. 

Finally, let us assume $H$ is almost simple and not of type (a) or (c). By \cite[Lemma 5.7]{AB}, the only possibility with $|H_0|^3 > |G_0|$ is the case $G = F_4(3)$ and $H = {}^3D_4(2).3$ as in part (v). By arguing as in the proof of the lemma, one checks that there are no additional subgroups with $|H_0|>q^{17}$.
\end{proof}

\begin{lem}
Suppose $G_0 = F_4(q)$, $q \geqs 3$ and $H$ is one of the subgroups in parts (ii) -- (v) of Lemma \ref{l:f4q17}. Then $b(G) \leqs 5$.
\end{lem}

\begin{proof}
If $H_0 = F_4(q^{1/2})$ or ${}^2F_4(q)$ then the result follows from the proof of \cite[Proposition 4.28]{BLS}. In each of the remaining cases, it is helpful to note that $|H_0|<q^{21}$ (see \cite[Table 5.1]{LSS} for the precise structure of $H_0$ in case (ii) of Lemma \ref{l:f4q17}).  Therefore, by arguing as in the proof of Lemma \ref{l:f4_bd}, it follows that the contribution from the elements $x \in G_0$ with $\dim x^{\G} \geqs 28$ is less than $q^{28}(q^{21}/q^{28})^5 = q^{-7}$. Similarly, the contribution from elements in the unipotent classes $A_1$ or $\tilde{A}_1$ (if $p=2$) is at most $q^{-2}$.

If $G = F_4(3)$ and $H = {}^3D_4(2).3$, it just remains to consider the contribution from the unipotent elements in the class $\tilde{A}_1$. Since $|x^G|>3^{21}$ and $i_3(H) < 3^{14}$, this is less than $3^{21}(3^{14}/3^{21})^5 = 3^{-14}$ and we conclude that $\what{Q}(G,5)<1$ as required. 

Next assume $H_0 = F_4(q_0)$, where $q=q_0^3$. If $x \in G_0$ is a unipotent element in the class $\tilde{A}_1$ ($p \ne 2$) or $(\tilde{A}_1)_2$ ($p=2$), then $|x^G|>\frac{1}{3}q^{22}$ and there are fewer than $2q_0^{22} = 2q^{22/3}$ such elements in $H$. Therefore, the contribution from these elements is less than $\frac{1}{3}q^{22}(6q^{-44/3})^5<q^{-44}$. As in the proof of Lemma \ref{l:f4_bd}, the contribution from involutory field or graph-field automorphisms is less than $2q^{-4}$. Similarly, the contribution from field automorphisms of odd prime order is less than $\frac{1}{2}q^{34}(4\log_2q.q^{-13})^5<q^{-23}$ and it follows that
$$\what{Q}(G,5)< q^{-2}+2q^{-4}+q^{-7}+q^{-23}+q^{-44}< 1.$$

Finally, suppose $p=2$ and $H = N_G(\bar{H}_{\s})$ with $\bar{H}^0=C_2^2$, so $q \geqs 4$ and $H_0 = {\rm Sp}_{4}(q) \wr S_2$ or ${\rm Sp}_{4}(q^2).2$. As in the previous case, the contribution from elements in $G \setminus G_0$ is less than $2q^{-4}+q^{-23}$, so it just remains to consider the unipotent elements in the class $(\tilde{A}_1)_2$. Here $|x^G|>q^{22}$ and it is straightforward to show that $i_2(H_0) < 2q^{12}$. For example, if $H_0 =  {\rm Sp}_{4}(q^2).2$ then
\begin{align*}
i_2(H_0) & = |b_1^{{\rm Sp}_{4}(q^2)}| + |a_2^{{\rm Sp}_{4}(q^2)}| + |c_2^{{\rm Sp}_{4}(q^2)}| + 
|{\rm Sp}_{4}(q^2) : {\rm Sp}_{4}(q)| \\
& = (q^8-1)(q^4+1) + q^4(q^2+1)(q^4+1)
\end{align*}
and the claim follows (here we are using the notation from \cite{AS} for involutions in ${\rm Sp}_{4}(q^2)$). Therefore, the contribution to $\b$ from these elements is less than $q^{22}(2q^{-10})^5 < q^{-25}$ and the result follows.
\end{proof}

In order to complete the proof of Theorem \ref{t:main} when $G_0 = F_4(q)$ and $q \geqs 3$, we just need to consider the cases arising in part (i) of Lemma \ref{l:f4q17}. For $\bar{H} = A_1C_3$, the desired bound $b(G) \leqs 5$ is established in \cite[Lemma 4.27]{BLS}, so we may assume that $H = N_G(\bar{H}_{\s})$ with $\bar{H}^0 = B_4$, $C_4$ ($p=2$ only) or $D_4$.

\begin{lem}
If $G_0 = F_4(q)$, $q \geqs 3$ and $H = N_G(\bar{H}_{\s})$ with $\bar{H} = D_4.S_3$, then $b(G) \leqs 5$.
\end{lem}

\begin{proof}
As recorded in \cite[Table 5.1]{LSS}, there are two possibilities for $H_0$, namely 
$d^2.{\rm P\Omega}_{8}^{+}(q).S_3$ and ${}^3D_4(q).3$, where $d=(2,q-1)$. As usual, write $\what{Q}(G,5) = \a+\b+\gamma$.

We begin by estimating $\a$. Let $x \in G$ be a semisimple element of prime order $r$ and define $\bar{D}$, $z$ and $\delta(x)$ as before. As in \cite[(4.18)]{BLS}, we have
\begin{equation}\label{e:fr}
{\rm fpr}(x) < \frac{48(q+1)^z}{q^{\delta(x)+z-4}(q-1)^4}
\end{equation}

First assume $r=2$, so $\bar{D} = B_4$ or $A_1C_3$. In the first case, \cite[Theorem 2]{LLS2} gives ${\rm fpr}(x) \leqs 2q^{-5}$ and we deduce that the contribution to $\a$ is less than $2q^{16}(2q^{-5})^5 = 2^6q^{-9}$. Similarly, if $\bar{D} = A_1C_3$ then $\delta(x) \geqs 12$ by \cite[Theorem 2]{LLS} and thus \eqref{e:fr} gives ${\rm fpr}(x)<q^{-6}$, so the contribution from these involutions is at most $2q^{28}(q^{-6})^5 = 2q^{-2}$. 

Now assume $r$ is odd. If $\bar{D} = B_3T_1$ or $C_3T_1$ then \cite[Theorem 2]{LLS} gives $\delta(x) \geqs 12$ and we calculate that $G$ contains at most 
$$2|F_4(q):C_3(q)| = 2q^{15}(q^4+1)(q^{12}-1)$$ 
such elements. Therefore, the contribution to $\a$ is less than
$$2q^{15}(q^4+1)(q^{12}-1)\cdot \left(\frac{48(q+1)}{q^{9}(q-1)^4}\right)^5<q^{-2}.$$
If $\bar{D} = T_4$ then $\delta(x)=24$ and \eqref{e:fr} implies that ${\rm fpr}(x) < q^{-17}$, so the contribution from these elements is at most $q^{52}(q^{-17})^5 < q^{-33}$. 
For the remaining semisimple elements, we have $z \leqs 3$ and $\delta(x) \geqs 16$ by \cite[Theorem 2]{LLS}, so \eqref{e:fr} implies that ${\rm fpr}(x)<q^{-51/5}$ and it follows that the combined contribution to $\a$ from these elements is less than $q^{50}(q^{-51/5})^5 = q^{-1}$. We conclude that 
$$\a< q^{-1}+3q^{-2} + 2^6q^{-9} + q^{-33}.$$

Finally, in order to estimate $\b$ and $\gamma$ we can appeal directly to \cite[Lemma 4.26]{BLS}. Indeed, by inspecting the proof of this lemma, we immediately deduce that 
$$\b < 3q^{22}(2q^{-6})^5 + q^{48}(q^{-10})^5 < 2q^{-2}$$
and 
$$\gamma < 2q^{26}(q^{-6})^5 + \log_2q.q^{52}(6q^{-56/3})^5 < q^{-3}.$$
Therefore, 
$$\what{Q}(G,5) < q^{-1}+5q^{-2} + q^{-3}+2^6q^{-9} + q^{-33}$$
and the proof of the lemma is complete.
\end{proof}

\begin{lem}
If $G_0 = F_4(q)$, $q \geqs 3$ and $H = N_G(\bar{H}_{\s})$ with $\bar{H} = B_4$ or $C_4$, then $b(G) \leqs 5$.
\end{lem}

\begin{proof}
We proceed by estimating $\a$, $\b$ and $\gamma$ in turn. First let $x \in G$
be a semisimple element of prime order $r$ and define $\bar{D}$, $z$ and $\delta(x)$ in the usual manner. Note that $H_0 = 2.\O_{9}(q) = {\rm Spin}_{9}(q)$ if $q$ is odd, and $H_0 = {\rm Sp}_{8}(q)$ if $q$ is even (see \cite[Table 5.1]{LSS}).

First assume $r=2$, so $q$ is odd and $\bar{D} = B_4$ or $A_1C_3$. By inspecting \cite[Table 4.5.2]{GLS}, we see that $H_0$ contains three classes of involutions; the central involution $z$, together with involutions of type $z_1 = [-I_4,I_5]$ and $z_2 = [-I_8,I_1]$, with plus-type $(-1)$-eigenspaces. To determine the $G_0$-class of $z_1$ and $z_2$, it is sufficient to compute the dimension of the $1$-eigenspace of $z_i$ with respect to the adjoint action of $\bar{G}$ on the Lie algebra $V=\mathcal{L}(\bar{G})$ of $\bar{G}$. Indeed, as explained in \cite[Section 1.14]{Carter}, this coincides with the dimension of $C_{\bar{G}}(z_i)$, which uniquely determines the $G_0$-class of $z_i$. To do this, it is useful to observe that 
\begin{equation}\label{e:dec}
V{\downarrow}\bar{H} = \mathcal{L}(\bar{H}) \oplus U,
\end{equation} 
where $U$ is the $16$-dimensional spin module for $\bar{H}$. Furthermore, to determine the action of $z_i$ on $U$ we may assume that $z_i \in D_4<\bar{H}$ and we can use the fact that 
$U{\downarrow}D_4 = W_1 \oplus W_2$, where $W_1$ and $W_2$ are the $8$-dimensional spin modules for $D_4$. For example, $\dim C_{\bar{H}}(z_1) = 16$ and one checks that $z_1$ acts on $U$ as $[-I_8, I_8]$. Therefore, $\dim C_{V}(z_1) = 16+8$ and thus $C_{\bar{G}}(z_1) = A_1C_3$. In the same way, one checks that $C_{\bar{G}}(z_2) = B_4$. Therefore,
$$|z_1^G \cap H| = \frac{|{\rm SO}_{9}(q)|}{|{\rm SO}_{5}(q)||{\rm SO}_{4}^{+}(q)|} < 2q^{20},\;\; |z_1^G| >q^{28},$$
$$|z_2^G \cap H| = 1+\frac{|{\rm SO}_{9}(q)|}{|{\rm SO}_{8}^{+}(q)|} < 2q^{8},\;\; |z_2^G| >q^{16}$$
and we deduce that the contribution to $\a$ from involutions is less than $q^{-8}$.

Next assume $r=3$, so $q \geqs 4$ and $\bar{D} = B_3T_1$, $C_3T_1$ or $A_2\tilde{A}_2$ (see \cite[Table 4.7.1]{GLS}). As above, we can use \eqref{e:dec} to determine $\dim C_{\bar{G}}(z)$ for each element $z \in \bar{H}$ of order $3$. For example, suppose $z = [I_2, \omega I_3, \omega^2 I_3] \in D_4 < \bar{H}$, where $\omega \in K$ is a primitive cube root of unity. Here $\dim C_{\bar{H}}(z) = 12$ and we calculate that $z$ has Jordan form $[I_4,\omega I_6, \omega^2 I_6]$ on $U$. Therefore, $\dim C_V(z) = 12+4$ and thus $C_{\bar{G}}(z) = A_2\tilde{A}_2$. Similarly, one checks that if $z \in \bar{H}$ belongs to one of the other $\bar{H}$-classes, then $C_{\bar{G}}(z) = B_3T_1$ or $C_3T_1$. We deduce that there are less than $2q^{22}=a_1$ elements $x \in H_0$ of order $3$ with $\bar{D} = B_3T_1$ or $C_3T_1$. Moreover, $|x^G|>\frac{1}{2}(q+1)^{-1}q^{31} = b_1$. Similarly, there are fewer than $2q^{24}=a_2$ elements with $\bar{D} = A_2\tilde{A}_2$, and we have $|x^G|>\frac{1}{2}q^{36} = b_2$ in this case. Therefore, the contribution to $\a$ from elements of order $3$ is less than 
$\sum_{i=1}^{2}b_i(a_i/b_i)^5 <q^{-4}$.

For the remainder, we may assume that $r \geqs 5$. Note that $r$ divides $|Z(C_{G_0}(x))|$. In view of \cite[14.1]{GL}, the possibilities for $\bar{D}$ are as follows:
$$\begin{array}{l|cccccccc}
\bar{D} & B_3T_1 & C_3T_1 & A_2A_1T_1 & C_2T_2 & A_2T_2 & A_1^2T_2 & A_1T_3 & T_4 \\ \hline
\dim x^{\bar{G}} & 30 & 30 & 40 & 40 & 42 & 44 & 46 & 48 
\end{array}$$
By applying \cite[Lemma 4.5]{LLS2}, we deduce that 
\begin{equation}\label{e:ff}
{\rm fpr}(x) < \frac{6(q+1)^z}{q^{\delta(x)+z-4}(q-1)^4}
\end{equation}
(this bound is stated as (4.17) in \cite{BLS}).

First assume $\bar{D} = B_3T_1$ or $C_3T_1$, so $z=1$ and $\delta(x) \geqs 8$ by \cite[Theorem 2]{LLS}. As noted in the proof of the previous lemma, there are fewer than 
$2q^{15}(q^4+1)(q^{12}-1)$ such elements in $G$. Also observe that $|Z(C_{G_0}(x))| = q \pm 1$, so either $(r,q)=(5,4)$ or $q \geqs 8$, since $r \geqs 5$. Suppose $(r,q) = (5,4)$. By considering the proof of \cite[Lemma 4.5]{LLS2}, noting that $|(y^{\bar{H}})_{\s}|<q^{\dim y^{\bar{H}}}$ for all $y \in x^{\bar{G}}\cap \bar{H}$, we deduce that 
${\rm fpr}(x) < 6(q+1)q^{-9}$.
Therefore, the contribution to $\a$ is less than 
$$2q^{15}(q^4+1)(q^{12}-1)\left(6(q+1)q^{-9}\right)^5<q^{-1}.$$
By applying the upper bound in \eqref{e:ff}, one checks that the same conclusion holds if $q \geqs 8$.

The case $\bar{D} = A_2A_1T_1$ is similar. Here $z=1$, $q \geqs 4$ and $\delta(x) \geqs 12$ by \cite[Theorem 2]{LLS}. Furthermore, we calculate that there are fewer than $2q^{41}$ such elements in $G$. Therefore, by applying the bound in \eqref{e:ff}, we deduce that their contribution to $\a$ is less than 
$$2q^{41}\left(\frac{6(q+1)}{q^{9}(q-1)^4}\right)^5 < q^{-7}.$$

Next assume $\bar{D} = C_2T_2$ or $A_2T_2$, so $z=2$ and $\delta(x) \geqs 12$. First we calculate that there are less than $3q^{44}$ elements in $G$ of this form. By applying \eqref{e:ff}, we deduce that the contribution to $\a$ here is at most $q^{-2}$ if $q \geqs 4$. The case $q=3$ requires special attention. We claim that the contribution to $\a$ is still less than $q^{-2}$. To see this, first observe that $r \in \{5,7,13\}$, which implies that 
\begin{equation}\label{e:fbb}
|(y^{\bar{H}})_{\s}|< \left(\frac{q^3}{q^3-1}\right)q^{\dim y^{\bar{H}}}
\end{equation}
for all $y \in x^{\bar{G}}\cap \bar{H}$. Therefore, the contribution is less than 
$$3q^{44}\cdot \left(6\left(\frac{q^3}{q^3-1}\right)\left(\frac{q+1}{q}\right)^2q^{-12}\right)^5<q^{-2}$$
as claimed.

A similar argument applies if $\bar{D} = A_1^2T_2$. Here $G$ contains fewer than $3q^{46}$ such elements and as noted in the proof of \cite[Lemma 4.25]{BLS}, we have $\delta(x) \geqs 14$. By applying the bound in \eqref{e:ff}, it follows that the contribution to $\a$ is less than
$q^{-4}$.

To complete the analysis of semisimple elements, we may assume that $\bar{D} = A_1T_3$ or $T_4$. In the latter case, $\delta(x) = 16$ and \eqref{e:ff} implies that the contribution to $\a$ is less than 
$$q^{52}\left(\frac{6(q+1)^4}{q^{16}(q-1)^4}\right)^5 < q^{-7}.$$
Similarly, if $\bar{D} = A_1T_3$ then $\delta(x) \geqs 14$ and we calculate that there are fewer than $2q^{49}$ such elements in $G$. Using \eqref{e:ff}, we see that the contribution is at most $q^{-7}$ if $q \geqs 4$. Finally, if $q=3$ then $r \in \{5,7,13\}$ and \eqref{e:fbb} holds, hence 
$${\rm fpr}(x) < 6\left(\frac{q^3}{q^3-1}\right)\left(\frac{q+1}{q}\right)^3q^{-14}$$
and the contribution to $\a$ is less than $q^{-7}$.

Bringing together the above bounds, we conclude that
$$\a< q^{-1}+q^{-2}+2q^{-4}+3q^{-7}+q^{-8}.$$

Next let us consider $\b$. First observe that the $\bar{G}$-class of each unipotent element in $\bar{H}$ is determined in \cite[Section 4.4]{Lawunip}, which allows us to compute very accurate fixed point ratio estimates (this information is also presented in the proof of \cite[Lemma 4.6]{LLS} when $p=2$). For now, let us assume $p$ is odd and let $x \in G$ be a unipotent element of order $p$. As in the proof of \cite[Lemma 4.25]{BLS}, the contribution to $\b$ from the elements with $\dim x^{\bar{G}} \leqs 30$ is less than $\sum_{i=2}^{5}c_id_i^5<2q^{-2}$, where the $c_i$ and $d_i$ terms are recorded in \cite[Table 10]{BLS}. Now assume $\dim x^{\bar{G}}>30$, in which case the proof of \cite[Lemma 4.25]{BLS} gives ${\rm fpr}(x) < 3q^{-10}$. By considering the $\bar{G}$-class of each unipotent element in $H$, together with the class sizes in \cite[Table 22.2.4]{LSbook}, we calculate that there are fewer than $q^{42}$ unipotent elements $x \in G$ with ${\rm fpr}(x)>0$ that are not in the $\bar{G}$-class $F_4(a_1)$. Therefore, the contribution from these elements is less than $q^{42}(3q^{-10})^5 \leqs q^{-3}$. Finally, if $x$ is in $F_4(a_1)$ then $|x^G|>\frac{1}{4}q^{46}$ and there are fewer than $q^{32}$ such elements in $H$, hence the contribution from these elements is less than $\frac{1}{4}q^{46}(4q^{-14})^5<q^{-18}$. We conclude that 
$\b < 2q^{-2}+q^{-3}+q^{-18}$ when $p$ is odd.

Now assume $p=2$, so $q \geqs 4$. The $\bar{G}$-class of each involution in $\bar{H}$ is determined in the proof of \cite[Lemma 4.6]{LLS} and we deduce that the contribution to $\b$ from involutions in the classes labelled $\tilde{A}_1$, $(\tilde{A}_1)_2$ and $A_1\tilde{A}_1$ is less than  
$$\b < q^{16}(q^{-8})^5 + q^{22}(2q^{-6})^5 + q^{28}(2q^{-8})^5 < q^{-5}.$$
For example, if $x$ is in $(\tilde{A}_1)_2$ then $x^{\bar{G}} \cap \bar{H} = y_1^{\bar{H}} \cup y_2^{\bar{H}}$, where $y_1$ and $y_2$ are involutions of type $c_2$ and $a_4$, respectively, in terms of the notation in \cite{AS}. This gives  
$$|x^G \cap H| = \frac{|{\rm Sp}_{8}(q)|}{|{\rm Sp}_{4}(q)|q^{12}} + \frac{|{\rm Sp}_{8}(q)|}{|{\rm Sp}_{4}(q)|q^{10}} < 2q^{16},\;\; |x^G| = \frac{|F_4(q)|}{q^{20}|{\rm Sp}_{4}(q)|}>q^{22}$$
and thus the contribution from these elements is less than $q^{22}(2q^{-6})^5$. Finally, if $x \in G$ is a long root element then ${\rm fpr}(x) = (q^4-q^2+1)^{-1}$ and so the contribution here is 
$$(q^4+1)(q^{12}-1) \cdot (q^4-q^2+1)^{-5} < q^{-3}.$$
Therefore, $\b < q^{-3}+q^{-5}$ when $p=2$. 

Finally, let us consider $\gamma$. By arguing as in the proof of \cite[Lemma 4.25]{BLS}, we see that the contribution from involutory field and graph-field automorphisms is at most $2q^{26}(q^{-6})^5 = 2q^{-4}$. Similarly, if $x \in G$ is a field automorphism of odd prime order $r$, then ${\rm fpr}(x) < 4q^{-16(1-1/r)}$ and we deduce that
$$\gamma < 2q^{-4} + 4q^{104/3}(4q^{-32/3})^5 + \log_2q.q^{52}(4q^{-64/5})^5 < q^{-3}.$$

In conclusion, 
$$\what{Q}(G,5) <q^{-1}+3q^{-2}+2q^{-3} + 2q^{-4}+3q^{-7}+q^{-8}+q^{-18}$$
and this completes the proof of the lemma.
\end{proof}

Finally, we handle the case $q=2$. 

\begin{lem}\label{l:f42}
If $G_0 = F_4(2)$, then $b(G) \leqs 5$. 
\end{lem}

\begin{proof}
Here $G = F_4(2)$ or $F_4(2).2$. Note that the character table of $G$ is available in the Atlas \cite{Atlas}, and the maximal subgroups of $G$ are determined in \cite{NW}. As usual, set $H_0 = H \cap G_0$ and $\what{Q}(G,5) = \a+\b+\gamma$.

First assume $|H_0| \leqs 2^{17}$ and let $x \in G$ be an element of prime order. If $x$ is not a long (or short) root element, then $|x^G|>2^{22}$ and so the contribution to $\widehat{Q}(G,5)$ from these elements is less than $2^{22}(2^{17}/2^{22})^5 = 2^{-3}$. Now $G$ contains exactly $139230$ long and short root elements and \cite[Theorem 2]{LLS2} gives ${\rm fpr}(x) \leqs 1/13$, so the remaining contribution is at most $139230(1/13)^5<2^{-1}$ and the result follows. For the remainder, we may assume that $|H_0| > 2^{17}$. The possibilities for $H_0$ can be read off from \cite[Table 1]{NW}.

First let us consider the case where $G = F_4(2)$ and $H = {\rm Sp}_{8}(2)$. As noted in the proof of \cite[Lemma 4.25]{BLS}, the Web Atlas \cite{WebAt} provides a faithful permutation representation of $G$ of degree $69888$, corresponding to the action of $G$ on the set of cosets of $H$. With the aid of {\sc Magma} \cite{Magma}, it is easy to find three elements $x_i \in G$ such that 
$$H \cap H^{x_1} \cap H^{x_2} \cap H^{x_3}  = 1,$$
which implies that $b(G) \leqs 4$. Similarly, by taking the normalizer of a maximal subgroup $O_{8}^{+}(2) < {\rm Sp}_{8}(2)$, we quickly deduce that $b(G) = 3$ when $G = F_4(2)$ and $H = \O_{8}^{+}(2).S_3$. The case $G = F_4(2)$ and $H = {}^3D_4(2).3$ can also be handled using {\sc Magma}. Here we construct ${}^3D_4(2)<G$ by finding a pair of generators $u ,v$ such that $|u|=2$, $|v|=9$, $|uv|=13$ and $|uv^2|=8$ (see the Web Atlas \cite{WebAt}). Then $N_G(\la u,v \ra) = {}^3D_4(2).3$ and it is easy to check that $b(G)=3$.

The cases where $G = F_4(2).2$ and $H = {\rm Sp}_{4}(4).4$ or $({\rm Sp}_{4}(2) \wr S_2).2$ can be handled in a similar fashion. Here $H_0$ is a maximal subgroup of ${\rm Sp}_{8}(2) < F_4(2)$, so we can use the previous representation to construct $G_0$ and $H_0$ as permutation groups of degree $69888$. It is now straightforward to determine the $G_0$-class of each prime order element in $H_0$, which allows us to calculate $\a$ and $\b$ precisely. In both cases, we find that $\a+\b < 2^{-32}$. For $H = {\rm Sp}_{4}(4).4$ we have $\gamma = 0$ (there are no involutions in $H \setminus H_0$) and the result follows. Now assume $H = ({\rm Sp}_{4}(2) \wr S_2).2$ and let $x \in G$ be an involutory graph automorphism. Here $|x^G|>2^{26}$ and thus $\gamma< 2^{26}(|H_0|/2^{26})^5 < 2^{-4}$. We conclude that $\what{Q}(G,5)<1$, as required.

It is also convenient to use {\sc Magma} when $H_0 = {}^2F_4(2)$, which is the centralizer in $G_0$ of an involutory graph-field automorphism. First we construct ${}^2F_4(2)' = \la u,v \ra$ as a subgroup of $G_0$, where $|u|=2$, $|v|=3$ and $|uv|=13$, and we take its normalizer to get ${}^2F_4(2)$. As before, it is now easy to compute $\a$ and $\b$ precisely, and one checks that $\a+\b < 2^{-33}$. Finally, if $G = F_4(2).2$ and $x \in G$ is an involutory graph-field automorphism, then we may assume that 
$H = \la x \ra \times {}^2F_4(2)$. Therefore, $|x^G \cap H| = 1+i_2({}^2F_4(2)) = 13456$ and it follows that $\gamma<2^{-37}$.

To complete the proof of the lemma, we may assume that $H_0 = {\rm L}_{4}(3).2$. Let $x \in H$ be an element of prime order $r$ and note that $r \in \{2,3,5,13\}$. If $r=3$ then $i_r(H) = 82160$ and $|x^G|>2^{29}$, so the contribution to $\a$ from these elements is less than $2^{29}(82160/2^{29})^5<2^{-34}$. In fact, by computing $i_5(H)$ and $i_{13}(H)$, we deduce that $\a<2^{-34}$. Now assume $r=2$. If $x$ is a long (or short) root element then \cite[Theorem 2]{LLS2} gives ${\rm fpr}(x) \leqs 2^{-4}$ and there are fewer than $2^{18}$ of these elements in $G$, so the contribution to $\b$ is less than $2^{18}(2^{-4})^5 = 2^{-2}$. If $x \in G$ is any other involution, then $|x^G|>2^{22}$ and we note that $i_2(H) \leqs i_2({\rm Aut}({\rm L}_{4}(3))) = 27639$. Therefore, $\b+\gamma < 2^{-2}+2^{-14}$ and the result follows.
\end{proof}

\vs

This completes the proof of Theorem \ref{t:main}. 

\subsection{Proof of Corollary \ref{t:main2}}\label{ss:c2} 

Let $G \leqs {\rm Sym}(\O)$ be a finite almost simple primitive group, with point stabilizer $H$ and socle $G_0$, which is a simple exceptional group of Lie type over $\mathbb{F}_q$. In view of Remark \ref{r:par}, we may assume that $H$ is non-parabolic. Moreover, we may assume that $G_0 \in \{E_7(q), E_6^{\e}(q), F_4(q)\}$ by \cite[Theorem 4]{BLS}. To complete the argument, we now combine the proofs of \cite[Theorem 4]{BLS} and Theorem \ref{t:main}. In every case, there is an explicit function $f(q)$ such that $\what{Q}(G,5) < f(q)$ and $f(q) \to 0$ as $q$ tends to infinity. The result follows.


\begin{thebibliography}{99}

\bibitem{AB} S.H. Alavi and T.C. Burness, \emph{Large subgroups of simple groups},  
J. Algebra \textbf{421} (2015), 187--233.

\bibitem{AS} M. Aschbacher and G.M. Seitz, \emph{Involutions in Chevalley groups over fields of even order}, Nagoya Math. J. \textbf{63} (1976), 1--91.

\bibitem{Babai} L. Babai, \emph{On the order of uniprimitive permutation groups}, Annals of Math. \textbf{113} (1981), 553--568. 

\bibitem{Bai} A.A. Baikalov, \emph{Intersection of conjugate solvable subgroups in symmetric groups}, Algebra and Logic \textbf{56} (2017), 87--97.
 
\bibitem{BCam} R.F. Bailey and P.J. Cameron, \emph{Base size, metric dimension and other invariants of groups and graphs}, Bull. London Math. Soc. \textbf{43} (2011), 209--242.

\bibitem{BBR} J. Ballantyne, C. Bates and P. Rowley, \emph{The maximal subgroups of $E_7(2)$}, LMS J. Comput. Math. \textbf{18} (2015), 323--371.

\bibitem{Bochert} A. Bochert, \emph{{\"{U}}ber die {Z}ahl verschiedener {W}erte, die eine {F}unktion
  gegebener {B}uchstaben durch {V}ertauschung derselben erlangen kann}, Math. Ann. \textbf{33} (1889), 584--590.

\bibitem{Magma} W. Bosma, J. Cannon and C. Playoust,
\emph{The {\sc Magma} algebra system I: The user language},
J. Symbolic Comput. \textbf{24} (1997), 235--265.

\bibitem{Bou} N. Bourbaki, \emph{Groupes et Algebr\`{e}s de Lie (Chapters 4,5 and 6)},
Hermann, Paris (1968).

\bibitem{Bur7}
T.C. Burness, \emph{On base sizes for actions of finite classical groups}, J.
  London Math. Soc. \textbf{75} (2007), 545--562.

\bibitem{BGS}
T.C. Burness, R.M. Guralnick and J. Saxl, \emph{On base sizes for symmetric groups}, Bull. London Math. Soc. \textbf{43} (2011), 386--391.

\bibitem{BGS3}
T.C. Burness, R.M. Guralnick and J. Saxl, \emph{Base sizes for {$\mathcal{S}$}-actions of finite classical groups}, Israel J. Math. \textbf{199} (2014), 711--756.

\bibitem{BGS4}
T.C. Burness, R.M. Guralnick and J. Saxl, \emph{On base sizes for algebraic groups}, J. Eur. Math. Soc. (JEMS) \textbf{19} (2017), 2269--2341.

\bibitem{BGS2}
T.C. Burness, R.M. Guralnick and J. Saxl, \emph{Base sizes for geometric actions of finite classical groups}, in preparation.

\bibitem{BLS} T.C. Burness, M.W. Liebeck and A. Shalev, \emph{Base sizes for simple groups and a conjecture of Cameron}, Proc. London Math. Soc. \textbf{98} (2009), 116--162.

\bibitem{BOB} T.C. Burness and E.A. O'Brien and R.A. Wilson, \emph{Base sizes for sporadic simple groups}, Israel J. Math. \textbf{177} (2010), 307--334.

\bibitem{CamPG} P.J. Cameron, \emph{Permutation Groups}, London Math. Soc. Student Texts \textbf{45}, Cambridge University Press, 1999.

\bibitem{CK} 
P.J. Cameron and W.M. Kantor, \emph{Random permutations: some group-theoretic aspects}, Combin. Probab. Comput. \textbf{2} (1993), 257--262.

\bibitem{Carter} R.W. Carter, \emph{Finite Groups of Lie Type: Conjugacy Classes and Complex Characters}, John Wiley, London, 1985.

\bibitem{CLSS}
A.M. Cohen, M.W. Liebeck, J. Saxl and G.M. Seitz, \emph{The local maximal subgroups of exceptional groups of {L}ie type}, Proc. London Math. Soc. \textbf{64} (1992), 21--48.

\bibitem{CW} A.M. Cohen and D.B. Wales, \emph{Finite subgroups of $F_4(\mathbb{C})$ and $E_6(\mathbb{C})$}, Proc. London Math. Soc. \textbf{74} (1997), 105--150.

\bibitem{Atlas} J.H. Conway, R.T. Curtis, S.P. Norton, R.A. Parker and R.A. Wilson, \emph{Atlas of {F}inite {G}roups}, Oxford University Press, 1985.

\bibitem{FJ} P. Fleischmann and I. Janiszczak, \emph{The semisimple conjugacy classes of finite groups of Lie type $E_6$ and $E_7$}, Comm. Algebra \textbf{21} (1993), 93--161.

\bibitem{GL} D. Gorenstein and R. Lyons, \emph{The local structure of finite groups of characteristic $2$ type}, Mem. Amer. Math. Soc. \textbf{42}, 1983.

\bibitem{GLS}
D. Gorenstein, R. Lyons  and R. Solomon, \emph{The classification of the
  finite simple groups, Number 3}, Mathematical Surveys and Monographs,
  vol. 40, Amer. Math. Soc., 1998.

\bibitem{KL} P.B. Kleidman and M.W. Liebeck, \emph{The {S}ubgroup {S}tructure of the {F}inite {C}lassical {G}roups}, London Math. Soc. Lecture Note Series, vol. 129, Cambridge University Press, 1990.

\bibitem{KW} P.B. Kleidman and R.A. Wilson, \emph{The maximal subgroups of $E_6(2)$ and 
${\rm Aut}(E_6(2))$}, Proc. London Math. Soc. \textbf{60} (1990), 266--294.

\bibitem{Law} R. Lawther, \emph{Jordan block sizes of unipotent elements in exceptional algebraic groups}, Comm. Algebra \textbf{23} (1995), 4125--4156.

\bibitem{Lawunip} R. Lawther, \emph{Unipotent classes in maximal subgroups of exceptional algebraic groups}, J. Algebra \textbf{322} (2009), 270--293.

\bibitem{LLS} R. Lawther, M.W. Liebeck and G.M. Seitz, \emph{Fixed point spaces in actions of exceptional algebraic groups}, Pacific J. Math. \textbf{205} (2002), 339--391.

\bibitem{LLS2} R. Lawther, M.W. Liebeck and G.M. Seitz, \emph{Fixed point ratios in actions of finite exceptional groups of Lie type}, Pacific J. Math. \textbf{205} (2002), 393--464.

\bibitem{L10}
M.W. Liebeck, \emph{On minimal degrees and base sizes of primitive permutation groups}, Arch. Math. \textbf{43} (1984), 11--15.

\bibitem{LSS} M.W. Liebeck, J. Saxl and G.M. Seitz, \emph{Subgroups of maximal rank in finite exceptional groups of Lie type}, Proc. London Math. Soc. \textbf{65} (1992), 297--325.

\bibitem{LS90}
M.W. Liebeck and G.M. Seitz, \emph{Maximal subgroups of exceptional groups of Lie type, finite and algebraic}, Geom. Dedicata \textbf{36} (1990), 353--387.

\bibitem{LSbook} M.W. Liebeck and G.M. Seitz, \emph{Unipotent and {N}ilpotent {C}lasses in {S}imple {A}lgebraic {G}roups and {L}ie {A}lgebras}, Amer. Math. Soc. Monographs and Surveys series, volume 180, 2012.

\bibitem{LSh2} M.W. Liebeck and A. Shalev, \emph{Simple groups, permutation groups, and probability}, J. Amer. Math. Soc. \textbf{12} (1999), 497--520.

\bibitem{Lu2} F. L\"{u}beck, \emph{Centralizers and numbers of semisimple classes in exceptional groups of Lie type}, \\
\texttt{http://www.math.rwth-aachen.de/$\sim$Frank.Luebeck/chev/CentSSClasses}

\bibitem{Lub} F. L\"{u}beck, \emph{Generic {C}omputations in {F}inite {G}roups of {L}ie {T}ype}, book in preparation.

\bibitem{Lus} G. Lusztig, \emph{Character sheaves V}, Adv. Math. \textbf{61} (1986), 103--155.

\bibitem{Maroti} A. Mar\'{o}ti, \emph{On the order of primitive groups}, J. Algebra \textbf{258} (2002), 631--640.

\bibitem{Kou} V.D. Mazurov and E.I. Khukhro, \emph{Unsolved problems in group theory: The Kourovka notebook, no. 18 (English version)}, 2017 (arxiv:1401.0300).

\bibitem{NNOW} M. Neunh\"{o}ffer, F. Noeske, E.A. O'Brien and R.A. Wilson, \emph{Orbit invariants and an application to the Baby Monster}, J. Algebra \textbf{341} (2011), 297--305.

\bibitem{NW} S.P. Norton and R.A. Wilson, \emph{The maximal subgroups of $F_4(2)$ and its automorphism group}, Comm. Algebra \textbf{17} (1989), 2809--2824.
 
\bibitem{PS} C.E. Praeger and J. Saxl, \emph{On the orders of primitive permutation groups},  Bull. London Math. Soc. \textbf{12} (1980), 303--307.
 
\bibitem{Seress_book}
\'{A}. Seress, \emph{Permutation Group Algorithms}, Cambridge Tracts in Math. \textbf{152}, Cambridge University Press, 2003.
 
\bibitem{Shinoda} K. Shinoda, \emph{The conjugacy classes of Chevalley groups of type $(F_4)$ over finite fields of characteristic $2$}, J. Fac. Sci. Univ. Tokyo \textbf{21} (1974), 133--159.

\bibitem{Shoji} T. Shoji, \emph{The conjugacy classes of Chevalley groups of type $(F_4)$ over finite fields of characteristic $p \ne 2$}, J. Fac. Sci. Univ. Tokyo \textbf{21} (1974), 1--17.

\bibitem{Sims} C.C. Sims, \emph{Computation with permutation groups}, Proc. Second Sympos. on Symbolic and Algebraic Manipulation, (ACM, New York), pp.23--28, 1971.

\bibitem{Vdo} E.P. Vdovin, \emph{On the base size of a transitive group with solvable point stabilizer}, J. Algebra Appl. \textbf{11} (2012), 1250015, 14 pp.

\bibitem{Vdovin} E.P. Vdovin, \emph{On intersections of solvable Hall subgroups in finite simple exceptional groups of Lie type}, Tr. Inst. Mat. Mekh. \textbf{19} (2013), 62--70; translation in Proc. Steklov Inst. Math. \textbf{285} (2014), suppl. 1, S183--S190.
 
\bibitem{W18} R.A. Wilson, \emph{Maximal subgroups of ${}^2E_6(2)$ and its automorphism groups}, preprint, 2018 (arxiv:1801.08374).

\bibitem{WebAt}
R.A. Wilson et~al., \emph{A {W}orld-{W}ide-{W}eb {A}tlas of finite group
  representations}, \\
  \texttt{http://brauer.maths.qmul.ac.uk/Atlas/v3/}
\end{thebibliography}
\end{document}